# Relationships among Interpolation Bases of Wavelet Spaces and Approximation Spaces

Zhiguo Zhang

*College of Automation, University of Electronic Science and Technology of China, Chengdu, China*

Email: zhiguozhang@uestc.edu.cn

Address: School of Automation, University of Electronic Science and Technology of China, Xi Yuan Road, Chengdu, Sichuan 610054, China

Mark A. Kon

*Mathematics and Statistics, Boston University, Boston, USA*

# Relationships among Interpolation Bases of Wavelet Spaces and Approximation Spaces

*Zhiguo Zhang* □   *College of Automation, University of Electronic Science and Technology of China, Chengdu, China. Email: zhiguozhang@uestc.edu.cn*

*Mark A. Kon* □   *Mathematics and Statistics, Boston University, Boston, USA*

Abstract— A multiresolution analysis is a nested chain of related approximation spaces. This nesting in turn implies relationships among interpolation bases in the approximation spaces and their derived wavelet spaces. Using these relationships, a necessary and sufficient condition is given for existence of interpolation wavelets, via analysis of the corresponding scaling functions. It is also shown that any interpolation function for an approximation space plays the role of a special type of scaling function (an interpolation scaling function) when the corresponding family of approximation spaces forms a multiresolution analysis. Based on these interpolation scaling functions, a new algorithm is proposed for constructing corresponding interpolation wavelets (when they exist in a multiresolution analysis). In simulations, our theorems are tested for several typical wavelet spaces, demonstrating our theorems for existence of interpolation wavelets and for constructing them in a general multiresolution analysis.

## 1.    INTRODUCTION

Classic Shannon sampling theory is the foundation of discrete signal processing. However, since Shannon sampling often cannot deal with non-bandlimited signals in a satisfying way, reconstructing a non-bandlimited function from its samples has become an important topic both in the mathematical and signal processing literature.

Sampling in an approximation space of a multiresolution analysis is an example of recovery of non-bandlimited functions, also known as multiscale interpolation. Two important multiscale interpolation theories have been proposed respectively by D. L. Donoho [1] and G. G. Walter [2]. In the first theory (that of Donoho), interpolation wavelets are linked to autocorrelations of orthogonal wavelets, and approximation spaces are possibly non-orthogonal direct sums of wavelet spaces [1, 3-6]. In the second (that of Walter), interpolation wavelets are constructed from kernel functions of approximation spaces, and wavelet spaces are orthogonal to each other [2, 6-9]. In this paper, we restrict our attention to the second multiscale interpolation theory (that of Walter).

In [2], an orthonormal basis is associated to the existence of an interpolation basis for an approximation space $V_0$ in a multiresolution analysis. Using the discrete-time Fourier transform, Walter discusses invertibility of a sampling operator from an approximation space $V_0$ onto $l^2(\mathbb{R})$ and derives an algorithm for constructing an interpolation function from such an orthonormal basis [2]. Motivated by these results above, invertibility of the sampling operator and construction of an interpolation function were two important topics in a number of subsequent papers. In [6], a biorthogonal basis takes the place of an orthonormal basis in the construction of interpolation functions; more importantly, the Poisson summation formula is introduced for determining invertibility of a sampling operator. Based on the results of Walter and Aldoubi, Wen Chen [8, 9] suggests that a necessary and sufficient condition for an interpolation basis to exist in $V_0 \in L^2(\mathbb{R})$ is that the series $1/\sum_n \phi(n)e^{-inw}$ converges to an element of $L^2[0, 2\pi]$, where $\{\phi(x-n)\}_n$ is a Riesz basis of $V_0$.

In addition, irregular sampling theories are discussed via combinations of the Zak transform with wavelet sampling methods in [8, 9]. These irregular sampling theories focus on perturbations of regular sampling and use techniques similar to Kadec's 1/4 theorem [8, 9]. Since this, the sampling theories proposed by Walter and Chen have become an important basis for extensions of irregular sampling theories [10, 11], special sampling theories [12-16] and estimation of aliasing [17, 18] in approximation spaces.

In fact, the sampling theory of Walter and Chen goes further than just construction of an interpolation basis for an approximation space $V_0$. In [2, 8, 9], although this is not discussed directly, it is indicated that the invertibility of the sampling operator, together with form of the reproducing kernel

$$q(x,y) = \sum_{n \in \mathbb{Z}} \phi(x-n)\tilde{\phi}(y-n), \qquad (1)$$

plays a critical role in interpolation in $V_0$; here $\{\tilde{\phi}(x-n)\}_n$ denotes the dual basis of $\{\phi(x-n)\}_n$.

From the kernel theory, sampling $f_s(x) \in V_0$ at a point $m \in \mathbb{Z}$ can always be represented as $f_s(m) = \int_{-\infty}^{+\infty} q(x,m) f_s(x) dx$. The results in [2, 8, 9] are important in that they confirm that the sequence $\{q(x,m)\}_{m \in \mathbb{Z}}$ forms a Riesz basis for $V_0$ if and only if the sampling operator $T_{V_0}(f_s) = \{f_s(k)\}_{k \in \mathbb{Z}}$ is an invertible linear transformation from the approximation space $V_0$ to $l^2(\mathbb{R})$. Furthermore, the Riesz basis $\{\tilde{S}_m^\phi(x) = q(x,m)\}_{m \in \mathbb{Z}}$ consists of translates of a fixed function, i.e. $\{\tilde{S}_m^\phi(x) = \tilde{S}^\phi(x-m)\}_{m \in \mathbb{Z}}$. This Riesz basis is biorthogonal to an interpolation basis $\{S_n^\phi(x)\}_{n \in \mathbb{Z}}$ in $V_0$, so the invertibility of the sampling operator $T_{V_0}(f_s)$ becomes a necessary and sufficient condition for an interpolation basis $\{S^\phi(x-n)\}_{n\in\mathbb{Z}}$ to exist in $V_0$ [2, 8, 9]. In [2], there are examples involving several familiar multiresolution analyses where the sampling operator $T_{V_0}(f_s) = \{f_s(k)\}_{k\in\mathbb{Z}}$ is not an invertible transformation from $V_0$ onto $l^2(\mathbb{R})$ and there do not exist interpolation bases in $V_0$.

Obviously, the approaches in [2, 8, 9] are not only suitable for the approximation space $V_0$, but also for any reproducing kernel Hilbert space (RKHS) with a kernel given by $\sum_{n\in\mathbb{Z}}\gamma(x-n)\tilde{\gamma}(y-n)$. Here $\{\gamma(x-n)\}_n$ denotes a Riesz basis and $\{\tilde{\gamma}(x-n)\}_n$ is the dual basis of $\{\gamma(x-n)\}_n$. A wavelet space $W_0$ is also a RKHS if its reproducing kernel $\sum_{n\in\mathbb{Z}}\psi(x-n)\tilde{\psi}(y-n)$ converges, where $\psi(x)$ is a wavelet in $W_0$ and $\tilde{\psi}(x)$ is the dual wavelet of $\psi(x)$. Hence, we consider whether there is a relationship between interpolation bases in approximation spaces and in wavelet spaces. It is well known that a multiresolution analysis is a nested chain of approximation spaces. The relationships among successive approximation spaces indicate that the interpolation bases of approximation spaces should be closely related to those of wavelet spaces. However, these relationships among different spaces are seldom considered in wavelet sampling theories. Thus, from the viewpoint of multiresolution approximation, wavelet sampling theories are still far from achieving what the Shannon theorem has accomplished.

In this paper, we do not focus on construction of an interpolation basis in a single space, but on relationships among those such bases in different wavelet spaces and in approximation spaces. Let $\hat{\phi}(w)$ be the Fourier transform of $\phi(x)$. Under the assumption that the functions $\sum_{k\in\mathbb{Z}}\phi(k/2)e^{-iwk/2}$ and $\sum_{k\in\mathbb{Z}}\hat{\phi}(w+4k\pi)$ converge, we will give a necessary and sufficient condition for existence of an interpolation basis in a wavelet space via analysis of a scaling function $\phi(x)$. We will also show that the interpolation functions of approximation spaces are in fact special scaling functions, so there exists a pair of reconstruction filters corresponding to interpolation scaling functions and interpolation wavelets. Thus interpolation bases in approximation spaces and in wavelet spaces are in fact related to each other. Finally, using this pair of reconstruction filters, an algorithm is derived for constructing an interpolation wavelet from an interpolation scaling function. Since it is relatively difficult to obtain analytical expressions for wavelets, our algorithm simplifies the construction of interpolation wavelets.

## 2. INTERPOLATION SCALING FUNCTIONS AND THEIR TWO-SCALE SYMBOLS

We will use a standard notation throughout, which is listed in Appendix G. Let $W_j$ denote a wavelet space and $V_j$ an approximation space such that $V_{j+1} = V_j \oplus W_j$ so that $\{V_j\}_{j \in \mathbb{Z}}$ forms a multiresolution analysis (MRA) of $L^2(\mathbb{R})$ where $\oplus$ denotes an orthogonal sum. The function $\phi(x)$ is assumed to be a real scaling function for the approximation space $V_0$ and $\psi(x)$ is its corresponding wavelet.

In this paper, we discuss relationships among interpolation bases in $V_j$ and $W_j$ under the assumption that there exists an interpolation basis $\{S^\phi(x-k)\}_{k \in \mathbb{Z}}$ in the approximation space $V_0$. To simplify the discussion, we always suppose that interpolation points of $\{S^\phi(x-k)\}_{k \in \mathbb{Z}}$ are located at integers $k \in \mathbb{Z}$, so we have

$$f_s(x) = \sum_{k \in \mathbb{Z}} f_s(k) S^\phi(x-k), \tag{2}$$

for any element $f_s \in V_0$. However, the results of this paper can be easily extended to the condition that the interpolation points of $\{S^\phi(x-k)\}_{k \in \mathbb{Z}}$ are located at other positions by the Zak transform or coordinate translation.

In addition, we assume that the series $\sum_{k \in \mathbb{Z}} \phi(k/2) e^{-iwk/2}$ and $\sum_{k \in \mathbb{Z}} \hat{\phi}(w+4k\pi)$ converge pointwise everywhere so that the following identity

$$\frac{1}{2} \sum_{k \in \mathbb{Z}} \phi(k/2) e^{-iwk/2} = \sum_{k \in \mathbb{Z}} \hat{\phi}(w+4k\pi) \tag{3}$$

holds by the Poisson summation formula. Here and henceforth we normalize the Fourier transform as $\hat{\phi}(w) = \int_{-\infty}^{+\infty} \phi(x) e^{-iwx} dx$. Obviously, since $\{\phi(k)\}_{k \in \mathbb{Z}}$ and $\{\phi(k+1/2)\}_{k \in \mathbb{Z}}$ are subsequences of $\{\phi(k/2)\}_{k \in \mathbb{Z}}$ and since $\{e^{-iwk/2}\}_{k \in \mathbb{Z}}$ are linearly independent, it follows from (3) that the series $\sum_{k \in \mathbb{Z}} \phi(k) e^{-iwk}$ (or $\sum_{k \in \mathbb{Z}} \phi(k+1/2) e^{-iwk}$) and $\sum_{k \in \mathbb{Z}} \hat{\phi}(w+2k\pi)$ (or $\sum_{k \in \mathbb{Z}} (-1)^k \hat{\phi}(w+2k\pi)$) also converge pointwise everywhere (see (A2) and (A3) in Appendix A). Thus, by the Poisson summation formula, we also have

$$\begin{cases} \sum_{k \in \mathbb{Z}} \phi(k) e^{-iwk} = \sum_{k \in \mathbb{Z}} \hat{\phi}(w+2k\pi) \\ e^{-iw/2} \sum_{k \in \mathbb{Z}} \phi(k+1/2) e^{-iwk} = \sum_{k \in \mathbb{Z}} (-1)^k \hat{\phi}(w+2k\pi) \end{cases}. \tag{4}$$

In this section, we show that the above interpolation function $S^\phi(x)$ is a scaling function with a special two-scale symbol. In addition, we show that when the series $\sum_{k\in\mathbb{Z}}\phi(k/2)e^{-iwk/2}$ and $\sum_{k\in\mathbb{Z}}\hat{\phi}(w+4k\pi)$ converge, the series $\sum_{k\in\mathbb{Z}}S^\phi(k/2)e^{-iwk/2}$ and $\sum_{k\in\mathbb{Z}}\hat{S}^\phi(w+4k\pi)$ also converge. Thus, we can apply the Poisson summation formula to $\{S^\phi(k/2)\}_{k\in\mathbb{Z}}$ in our proofs.

Assume that $U_j(S^\phi)$ is the space spanned by $\{S^\phi(2^j x-k)\}_k$ with $j\in\mathbb{Z}$. From [6] (Section 7.6.1), a sequence $\{S^\phi(2^j x-k)\}_k$ forms an interpolation basis for $U_j(S^\phi)$ if and only if $\{S^\phi(x-k)\}_k$ is an interpolation basis for $U_0(S^\phi)$. Obviously, $U_0(S^\phi)$ is just $V_0$ and $U_j(S^\phi)$ is just $V_j$ in this paper. Hence, $S^\phi(x)$ is a scaling function.

On the other hand, for any element $f_s \in V_j$, $j\in\mathbb{Z}$, we have $f_s(x/2^j) \in V_0$ (Lemma 5.1 in [20]). Hence, it follows from (2) that

$$f_s(y/2^j) = \sum_{k\in\mathbb{Z}} f_s(k/2^j) S^\phi(y-k). \tag{5}$$

The change of variable $x = y/2^j$ in (5) gives

$$f_s(x) = \sum_{k\in\mathbb{Z}} f_s(k/2^j) S^\phi(2^j x - k). \tag{6}$$

Hence, it follows from (6) that $S^\phi(x)$ is an interpolation scaling function such that $\{S^\phi(2^j x - k)\}_{k\in\mathbb{Z}}$ have interpolation points located at $k/2^j$, $k\in\mathbb{Z}$.

Let $P_\phi(z)$ denote a two-scale symbol related to a scaling function $\phi(x)$. By definition (see equation (5.1.8) in [20]), a two-scale symbol is an $l^2$-sequence $\{p_k\}_{k\in\mathbb{Z}}$ such that

$$\phi(x) = \sum_{k\in\mathbb{Z}} p_k \phi(2x-k) \tag{7}$$

and hence

$$\hat{\phi}(w) = \frac{1}{2}\sum_{k\in\mathbb{Z}} p_k z^k \hat{\phi}(\frac{w}{2}) = P_\phi(z)\hat{\phi}(\frac{w}{2}), \tag{8}$$

where

$$P_\phi(z) = \frac{1}{2}\sum_{k\in\mathbb{Z}} p_k z^k, \ z = e^{-iw/2}. \tag{9}$$

In [2, 8, 9], it is shown that the Fourier transform $\hat{S}^\phi(w)$ of an interpolation function $S^\phi(x)$ is expressed in terms of the Fourier transform $\hat{\phi}(w)$ as

$$\hat{S}^{\phi}(w) = \frac{\hat{\phi}(w)}{\sum_{k \in \mathbb{Z}} \phi(k)e^{-iwk}} = \frac{\hat{\phi}(w)}{\sum_{k \in \mathbb{Z}} \hat{\phi}(w+2k\pi)} \qquad (10)$$

with

$$\frac{1}{\sum_{k \in \mathbb{Z}} \hat{\phi}(w+2k\pi)} \in L^2[0, 2\pi] \ .$$

By using (3) and (10), Proposition 1 below describes a relationship among the interpolation bases $\{S^{\phi}(2^j x - k)\}_{k \in \mathbb{Z}}$ of different approximation spaces $V_j$, $j \in \mathbb{Z}$. It states that there exists a special two-scale symbol $P_s(z)$ corresponding to the interpolation scaling function $S^{\phi}(x)$.

**Proposition 1.** Assume that $U_0(\gamma)$ is a reproducing kernel Hilbert space spanned by a Riesz basis $\{\gamma(x-k)\}_k$ and that $S^{\gamma}(x)$ is an interpolation function for $U_0(\gamma)$, whose Fourier transform $\hat{S}^{\gamma}(w) = \hat{\gamma}(w) / \sum_{k \in \mathbb{Z}} \gamma(k)e^{-iwk}$ with $\hat{\gamma}(w)$ the Fourier transform of $\gamma(x)$. If the series $\sum_{k \in \mathbb{Z}} \gamma(k/2)e^{-iwk/2}$ and $\sum_{k \in \mathbb{Z}} \hat{\gamma}(w+4k\pi)$ converge pointwise everywhere, then the series $\sum_{k \in \mathbb{Z}} S^{\gamma}(k/2)e^{-iwk/2}$ and $\sum_{k \in \mathbb{Z}} \hat{S}^{\gamma}(w+4k\pi)$ also do, so that we have

$$\frac{1}{2} \sum_{k \in \mathbb{Z}} S^{\gamma}(k/2) e^{-iwk/2} = \sum_{k \in \mathbb{Z}} \hat{S}^{\gamma}(w+4k\pi) \qquad (11)$$

by the Poisson summation formula.

Thus when there exists an interpolation basis $\{S^{\phi}(x-k)\}_k$ relative to the samples $\{f_s(k)\}_k$ (*interpolation points are located at $k \in \mathbb{Z}$*) in $V_0$, then

$$P_s(z) = \frac{1}{2} \sum_{k \in \mathbb{Z}} S^{\phi}(\frac{k}{2}) z^k = \sum_{k \in \mathbb{Z}} \hat{S}^{\phi}(w+4k\pi) \qquad (12)$$

is the two-scale symbol corresponding to the scaling function $S^{\phi}(x)$. Furthermore, we have

$$\hat{S}^{\phi}(w) = P_s(z) \hat{S}^{\phi}(w/2) = \sum_{k \in \mathbb{Z}} \hat{S}^{\phi}(w+4k\pi) \hat{S}^{\phi}(\frac{w}{2^n}) \qquad (13)$$

and

$$P_s(z) = \frac{P_{\phi}(z) \sum_{k \in \mathbb{Z}} \hat{\phi}(w/2 + 2k\pi)}{\sum_{n \in \mathbb{Z}} \hat{\phi}(w+2n\pi)} \ . \qquad (14)$$

**Proof.** See Appendix A for the proof.

Equation (12) indicates that the symbol $P_s(z)$ is thus generated by $P_s(z) = 1/2 \sum_{k \in \mathbb{Z}} S^{\phi}(k/2) z^k$, where the sequence $\{S^{\phi}(k/2)\}_k$ is obtained by sampling the interpolation function $S^{\phi}(2^j x)$ at intervals $1/2^{j+1}$. By Proposition 1, the interpolation function $S^\phi(x)$ is a scaling function corresponding to the two-scale symbol $P_s(z)$ in (12). Hence, $P_s(z)$ has connected the interpolation functions $S^\phi(2^j x)$ at all scales $1/2^j$.

Obvious as the result of Proposition 1 may seem, it is important for us to describe interpolation wavelets. Just from the result of Proposition 1, Section 3 indicates that an interpolation function of approximation space $V_j$ often corresponds to that of the wavelet space $W_j$.

**Remark.** If an interpolation function $S^\phi(x)$ is a scaling function, then a two-scale symbol $P_s(z)$ connects the interpolation functions $S^\phi(2^j x)$ at all resolutions $2^j$ by repeated application of (13). Hence, equation (13) implies that the interpolation bases $\{S^\phi(2^j x - k)\}_{k \in \mathbb{Z}}$ for all $V_j$ have a close relationship in an MRA. More importantly, it is shown in Proposition 1 that $P_s(z)$ can be presented as in terms of $\hat{S}^\phi(w)$ in (12), and is connected to $P_\phi(z)$ in (14). The following sections show that (12) and (14) are important for connecting existence of an interpolation wavelet $S^\psi(2^j x)$ to that of $S^\phi(x)$.

## 3. A SUFFICIENT CONDITION FOR EXISTENCE OF INTERPOLATION BASES IN WAVELET SPACES

### 3.1 Invertibility of the Sampling Operator from Wavelet Spaces to $l^2(\mathbb{R})$

In Section 2, we have shown that $S^\phi(x)$ is an interpolation scaling function, so there exists the two-scale symbol $P_s(z)$ corresponding to $S^\phi(x)$ as in (12).

Let

$$PE_s(w) = \overline{P_s(-z)}E_s(-z) + \overline{P_s(z)}E_s(z) \tag{15}$$

where

$$E_s(z) = E_s(e^{-iw/2}) = \sum_{k=-\infty}^{+\infty} \left| \hat{S}^\phi\left(\frac{w}{2} + 2k\pi\right) \right|^2. \tag{16}$$

This section shows that $PE_s(w)$ in (15) plays a critical role in existence of an interpolation basis $\{S^\psi(2^j x - k)\}_{k \in \mathbb{Z}}$ in $W_j$ such that

$$f_s(x) = \sum_{k=-\infty}^{+\infty} f_s(k/2^j + 1/2^{j+1}) S^\psi(2^j x - k) \tag{17}$$

for any $f_s \in W_j$ and

$$\int_{-\infty}^{+\infty} S^\psi(2^{j_1} x - k) S^\phi(2^{j_2} x - l) dx = 0, \text{ for } k, l, j_1, j_2 \in \mathbb{Z} \tag{18}$$

when there exists an interpolation basis $\{S^\phi(x-k)\}_{k \in \mathbb{Z}}$ relative to the samples $\{f_s(k)\}_{k \in \mathbb{Z}}$ in $V_0$.

It has been observed [2, 8, 9] that when $V_0$ is a reproducing kernel Hilbert space with its kernel given by (1) (i.e. when the series $\sum_{n \in \mathbb{Z}} \phi(x-n)\tilde{\phi}(y-n)$ converges) and when $T_{V_0} f_s = \{f_s(k)\}_k$ is an invertible transformation from $V_0$ onto $l^2(\mathbb{R})$, then $V_0$ has an interpolation basis of the form $\{S^\phi(x-k)\}_k$. This result can be generalized easily to any functional Hilbert space $U_j(\gamma)$ spanned by a Riesz basis $\{\gamma(2^j x - k)\}_{k \in \mathbb{Z}}$. This is summarized in

**Proposition 2.** Assume $U_j(\gamma)$ is a reproducing kernel Hilbert space with a kernel

$$\begin{aligned} q_j(x,y) &= 2^j q(2^j x, 2^j y) \\ &= 2^j \sum_{n=-\infty}^{+\infty} \gamma(2^j x - n)\tilde{\gamma}(2^j y - n) \end{aligned} \tag{19}$$

where $\{\tilde{\gamma}(2^j x - k)\}_{k \in \mathbb{Z}}$ is the dual basis of the Riesz basis $\{\gamma(2^j x - k)\}_{k \in \mathbb{Z}}$. Let a sampling operator $T_{U_j}$ associate $f_s(x) \in U_j(\gamma)$ to its samples located at $x = k/2^j$ with $k \in \mathbb{Z}$, i.e. $T_{U_j}(f_s) = \{f_s(k/2^j)\}_k$. Then if $T_{U_j}$ is an invertible linear transformation of $U_j(\gamma)$ onto $l^2(\mathbb{R})$, there exists an interpolation basis $\{S^\gamma(2^j x - k)\}_k$ in $U_j(\gamma)$. Furthermore, the inverse of $T_{U_j}$ is represented by

$$f_s(x) = \sum_{k \in \mathbb{Z}} f_s(k/2^j) S^\gamma(2^j x - k). \tag{20}$$

**Proof.** Since the proof is very similar to that of [2], [8] and [9], we will not elaborate on it.

Obviously, every $W_j$ is just a reproducing kernel space whose reproducing kernel has the form (19), so it follows from Proposition 2 that invertibility of sampling operator is critical for existence of an interpolation basis in a wavelet space.

Let

$$\tilde{T}_{W_{j-1}}(f_s) = \{f_s(k/2^{j-1} + 1/2^j)\}_k \tag{21}$$

denote a sampling operator from $W_{j-1}$ onto $l^2(\mathbb{R})$. In Theorem 1, we will show that $\tilde{T}_{W_{j-1}}$ is an invertible operator when $0 < |PE_s(w)|^2 < +\infty$. In order to prove this, we first discuss some properties of the dual basis $\{\tilde{S}^\phi(x-k)\}_{k\in\mathbb{Z}}$ of $\{S^\phi(x-k)\}_{k\in\mathbb{Z}}$.

Since $S^\phi(x)$ is a scaling function by Proposition 1 and since $\{\tilde{S}^\phi(x-k)\}_{k\in\mathbb{Z}}$ is the dual basis of $\{S^\phi(x-k)\}_{k\in\mathbb{Z}}$ in $V_0$, $\tilde{S}^\phi(x)$ is also a scaling function (see Theorems 3.25 and 5.22 in [20]). Hence, $\tilde{S}^\phi(x)$ has a two-scale symbol $P_{\tilde{S}^\phi}(z)$ as in (7) such that

$$\begin{cases} \tilde{S}^\phi(x) = \sum_{k\in\mathbb{Z}} \tilde{p}_k^s \tilde{S}^\phi(2x-k) \\ \hat{\tilde{S}}^\phi(w) = \frac{1}{2}\sum_{k\in\mathbb{Z}} \tilde{p}_k^s z^k \hat{\tilde{S}}^\phi(\frac{w}{2}) \end{cases} \quad (22)$$

where $\hat{\tilde{S}}^\phi(w)$ is the Fourier transform of $\tilde{S}^\phi(x)$. From (22), we have

**Lemma 1.** Assume that there exists an interpolation basis $\{S^\phi(x-k)\}_{k\in\mathbb{Z}}$ relative to the samples $\{f_s(k)\}_{k\in\mathbb{Z}}$ in $V_0$ and that there exist two constants $A_s$ and $B_s$ with $0 < A_s \leq B_s < +\infty$ such that

$$A_s \leq |PE_s(w)|^2 = \left|\overline{P_s(-z)}E_s(-z) + \overline{P_s(z)}E_s(z)\right|^2 \leq B_s. \quad (23)$$

Then, for any given set of coefficients $\{c_k^0\}_{k\in\mathbb{Z}} \in l^2(\mathbb{R})$, there exists two unique sets of coefficients $\{b_k^0\}_{k\in\mathbb{Z}}$, $\{d_k^0\}_{k\in\mathbb{Z}} \in l^2(\mathbb{R})$ such that

$$\sum_{k\in\mathbb{Z}} b_k^0 \tilde{S}^\phi(2^j x - k) = \sum_{k\in\mathbb{Z}} c_k^0 \tilde{S}^\phi(2^{j+1} x - 2k) + \sum_{k\in\mathbb{Z}} d_k^0 \tilde{S}^\phi(2^{j+1} x - 2k - 1), \quad (24)$$

where $\{\tilde{S}^\phi(x-k)\}_{k\in\mathbb{Z}}$ is the dual basis of $\{S^\phi(x-k)\}_{k\in\mathbb{Z}}$.

**Proof.** See Appendix B for the proof.

Obviously, on the right-hand side of (24), $\sum_{k\in\mathbb{Z}} c_k^0 \tilde{S}^\phi(2^{j+1}x-2k) + \sum_{k\in\mathbb{Z}} d_k^0 \tilde{S}^\phi(2^{j+1}x-2k-1)$ represents an element of $V_{j+1}$, while on the left-hand side of (24), $\sum_{k\in\mathbb{Z}} b_k^0 \tilde{S}^\phi(2^j x - k)$ represents an element of $V_j$. Hence, Lemma 1 indicates that the coefficients $\{c_k^0\}_{k\in\mathbb{Z}}$ in (24) determine uniquely an element of $V_j$. In other words, there exists a bijection between $\{c_k^0\}_{k\in\mathbb{Z}}$ and $\{d_k^0\}_{k\in\mathbb{Z}}$ in (24) when (23) holds.

Based on Lemma 1, we show in Theorem 1 that the sampling operator $\tilde{T}_{W_{j-1}}(f_s)$ is an invertible transformation from $W_{j-1}$ onto $l^2(\mathbb{R})$ if (23) holds.

**Theorem 1.** Assume that there exists an interpolation basis $\{S^\phi(x-k)\}_{k\in\mathbb{Z}}$ relative to the samples $\{f_s(k)\}_{k\in\mathbb{Z}}$ in $V_0$ and that (23) holds. Then, for any given sequence $\{F_k\}_{k\in\mathbb{Z}} \in l^2(\mathbb{R})$, there exists a unique element $f_s(x) \in W_{j-1}$ such that

$$\tilde{T}_{W_{j-1}}(f_s) = \left\{f_s(\frac{2k+1}{2^j})\right\}_{k\in\mathbb{Z}} = \{F_k\}_{k\in\mathbb{Z}}, \; k \in \mathbb{Z}. \tag{25}$$

**Proof.** See Appendix C for the proof.

Since Theorem 1 indicates that every $\{f_s(k/2^{j-1}+1/2^j)\}_k \in l^2(\mathbb{R})$ corresponds to a unique element of $W_{j-1}$ when (23) holds, it follows that the sampling operator $\tilde{T}_{W_{j-1}}$ is a bijection between $W_{j-1}$ and $l^2(\mathbb{R})$, i.e. $\tilde{T}_{W_{j-1}}$ is an invertible transformation from $W_{j-1}$ onto $l^2(\mathbb{R})$.

## 3.2 *Inverse of the Sampling Operator from Wavelet Spaces onto $l^2(\mathbb{R})$*

Now we focus on how to construct the inverse of $\tilde{T}_{W_{j-1}}$ based on a wavelet $\psi(x)$ in $W_0$. In addition, we will show that when $\sum_{k\in\mathbb{Z}} \phi(k/2)e^{-iwk/2}$ and $\sum_{k\in\mathbb{Z}} \hat{\phi}(w+4k\pi)$ converge in (3), $\sum_{k\in\mathbb{Z}} S^\psi(k/2)e^{-iwk/2}$ and $\sum_{k\in\mathbb{Z}} \hat{S}^\psi(w+4k\pi)$ do also, so that we can apply the Poisson summation formula to $\{S^\psi(k/2)\}_{k\in\mathbb{Z}}$.

From Proposition 2, the inverse of the sampling operator $T_{U_j}$ can be constructed using (20) when $T_{U_j}$ is invertible. However, Proposition 2 is only suitable for the samples of $f_s(x)$ located at $x = k/2^j$, $k \in \mathbb{Z}$. Observe that the sampling operator $\tilde{T}_{W_{j-1}}(f_s) = \{f_s(k/2^{j-1}+1/2^j)\}_{k\in\mathbb{Z}}$ locates the samples of $f_s(x)$ at $x = (2k+1)/2^j$, so that Proposition 2 cannot be applied to construct its inverse directly. In order to construct the inverse, the space

$$\tilde{W}_{j-1} = \left\{f(x)\big| f(x-1/2^j) \in W_{j-1}\right\} \tag{26}$$

is introduced in the following proposition, which involves an indirect application of Proposition 2.

**Proposition 3.** Let the sampling operator $\tilde{T}_{W_{j-1}}(f_s) = \{f_s(k/2^{j-1}+1/2^j)\}_{k\in\mathbb{Z}}$ be a bijection of $W_{j-1}$ onto $l^2(\mathbb{R})$ and $\{\psi(x-k)\}_k$ be a Riesz basis of $W_0$. Then,

(a) The inverse of the sampling operator $\tilde{T}_{W_{j-1}}$ is given by

$$f_s(x) = \sum_{k\in\mathbb{Z}} f_s(\frac{k}{2^{j-1}}+\frac{1}{2^j})S^\psi(2^{j-1}x-k), \tag{27}$$

where the Fourier transform of $S^{\psi}(x)$ is

$$\hat{S}^{\psi}(w) = \frac{\hat{\psi}(w)}{\sum_{k \in \mathbb{Z}} \psi(k+1/2)e^{-iwk}}. \tag{28}$$

(b)  Furthermore, let $Q_{\psi}(z) = 1/2 \sum_n q_n e^{-iwn/2}$ denote a two-scale symbol corresponding to the wavelet $\psi(x)$, i.e.

$$\begin{cases} \psi(x) = \sum_{k \in \mathbb{Z}} q_k \phi(2x-k) \\ \hat{\psi}(w) = Q_{\psi}(z)\hat{\phi}(w/2) \end{cases}, \tag{29}$$

where $\hat{\psi}(w)$ is the Fourier transform of $\psi(x)$. Then, if the series $Q_{\psi}(z)$ converge pointwise everywhere, it follows from (4) that $\sum_{k \in \mathbb{Z}} S^{\psi}(k/2)e^{-iwk/2}$ ( and $\sum_{k \in \mathbb{Z}} \psi(k/2)e^{-iwk/2}$ ) and $\sum_{k \in \mathbb{Z}} \hat{S}^{\psi}(w+4k\pi)$ (and $\sum_{k \in \mathbb{Z}} \hat{\psi}(w+4k\pi)$ ) also do. Thus

$$\begin{cases} \frac{1}{2}\sum_{k \in \mathbb{Z}} \psi(k/2)e^{-iwk/2} = \sum_{k \in \mathbb{Z}} \hat{\psi}(w+4k\pi) \\ \frac{1}{2}\sum_{k \in \mathbb{Z}} S^{\psi}(k/2)e^{-iwk/2} = \sum_{k \in \mathbb{Z}} \hat{S}^{\psi}(w+4k\pi) \end{cases} \tag{30}$$

holds by the Poisson summation formula. Hence, (28) can be rewritten as

$$\hat{S}^{\psi}(w) = \frac{e^{-iw/2}\hat{\psi}(w)}{\sum_k (-1)^k \hat{\psi}(w+2k\pi)}. \tag{31}$$

**Proof.** See Appendix D for the proof.

By statement (a), an interpolation basis $\{S^{\psi}(2^{j-1}x-k)\}_k$ of a wavelet space $W_j$ consists of translates of a fixed function $S^{\psi}(2^j x)$ and has the interpolation points located at $x = (2k+1)/2^j$.

By statement (b), under our assumption that $\sum_{k \in \mathbb{Z}} \phi(k/2)e^{-iwk/2}$ and $\sum_{k \in \mathbb{Z}} \hat{\phi}(w+4k\pi)$ in (3) and (4) converge pointwise everywhere, equation (30) holds for $\{\psi(k/2)\}_{k \in \mathbb{Z}}$ and $\{S^{\psi}(k/2)\}_{k \in \mathbb{Z}}$ when $Q_{\psi}(z)$ also converges pointwise everywhere. Hence, the Fourier transform of $S^{\psi}(x)$ in (28) can be written as in (31). In general, since $Q_{\psi}(z)$ always converges pointwise everywhere, we will directly use (30) and (31) in the following sections.

Theorem 1 and Proposition 3 indicate that when there exists the interpolation basis $\{S^{\phi}(x-k)\}_{k \in \mathbb{Z}}$ in $V_0$, equation (23) is a sufficient condition for the interpolation basis $\{S^{\psi}(2^j x-k)\}_{k \in \mathbb{Z}}$ to exist in $W_j$.

This implies that there is a close relationship between existence of interpolation bases in an approximation space $V_j$ and in the corresponding wavelet space $W_j$.

In the next section, we show that (23) is not only a sufficient condition, but also a necessary one for $\{S^\psi(2^j x - k)\}_{k \in \mathbb{Z}}$ to exist in $W_j$. Furthermore, an interpolation wavelet $S^\psi(x)$ is connected to the interpolation scaling function $S^\phi(x)$ by a special pair of construction filters, so that $S^\psi(x)$ can be constructed from $S^\phi(x)$ directly instead of a wavelet $\psi(x)$ as in (28).

**Remark.** Lemma 1 indicates that $\{c_k^0\}_{k \in \mathbb{Z}}$ determine $\{d_k^0\}_{k \in \mathbb{Z}}$ uniquely in (24) when (23) holds. This result forms an important condition for Theorem 1, which implies that the sampling operator $\tilde{T}_{W_j}$ is always invertible when $\{S^\phi(x - k)\}_{k \in \mathbb{Z}}$ exists in $V_0$. Hence, it follows from Propositions 2 and 3 that when $\{S^\phi(x - k)\}_{k \in \mathbb{Z}}$ exists in $V_0$, equation (23) is a sufficient condition for the existence of the interpolation basis $\{S^\psi(2^j x - k)\}_{k \in \mathbb{Z}}$ relative to samples $\{f_s(k/2^j + 1/2^{j+1})\}_{k \in \mathbb{Z}}$ in a wavelet space $W_j$. This implies that the existence of $\{S^\psi(2^j x - k)\}_{k \in \mathbb{Z}}$ has a close relation with that of $\{S^\phi(x - k)\}_{k \in \mathbb{Z}}$. In addition, it is shown in Proposition 3 that when $\sum_{k \in \mathbb{Z}} \phi(k/2) e^{-iwk/2}$ and $\sum_{k \in \mathbb{Z}} \hat{\phi}(w + 4k\pi)$ in (3) converge, the Fourier transform $\hat{S}^\psi(w)$ of $S^\psi(x)$ always can be represented as the form (31). It will be shown in the next section that (31) is important for us to construct an interpolation wavelet $S^\psi(x)$ from the interpolation scaling function $S^\phi(x)$.

## 4. A NECESSARY CONDITION FOR EXISTENCE OF INTERPOLATION BASES IN WAVELET SPACES

Section 2 states that $S^\phi(x)$ must be a scaling function with the corresponding two-scale symbol $P_s(z) = \sum_{k \in \mathbb{Z}} \hat{S}^\phi(w + 4k\pi)$ when the interpolation functions $\{S^\phi(x - k)\}_{k \in \mathbb{Z}}$ form a Riesz basis for $V_0$. In Section 3, Theorem 1 and Proposition 3 state that there exists the interpolation basis $\{S^\psi(2^j x - k)\}_k$ in every $W_j$, whenever there exists the interpolation basis $\{S^\phi(x - k)\}_{k \in \mathbb{Z}}$ in $V_0$ and (23) holds.

In this section, we show that if there exists the interpolation basis $\{S^\phi(x-k)\}_{k\in\mathbb{Z}}$ in $V_0$, then (23) is also a necessary condition for the existence of $\{S^\psi(2^j x-k)\}_k$ in $W_j$. Furthermore, we propose an algorithm to construct the interpolation function $S^\psi(x)$ from $S^\phi(x)$ directly.

Suppose that there exist the interpolation basis $\{S^\psi(2^j x-k)\}_{k\in\mathbb{Z}}$ relative to $\{f_s(k/2^j + 1/2^{j+1})\}_{k\in\mathbb{Z}}$ in $W_j$ and the interpolation basis $\{S^\phi(x-k)\}_{k\in\mathbb{Z}}$ relative to $\{f_s(k)\}_{k\in\mathbb{Z}}$ in $V_0$.

Since Proposition 1 indicates that $S^\phi(x)$ is a scaling function when $\{S^\phi(x-k)\}_{k\in\mathbb{Z}}$ forms an interpolation basis for $V_0$ and since $W_0 \subset V_1$, we have

$$S^\psi(x) = \sum_{k=-\infty}^{+\infty} S^\psi(k/2) S^\phi(2x-k), \tag{32}$$

with $\{S^\psi(k/2)\}_{k\in\mathbb{Z}} \in l^2(\mathbb{R})$.

By taking the Fourier transform on both sides in (32), we obtain

$$\hat{S}^\psi(w) = \frac{1}{2}\sum_{k\in\mathbb{Z}} S^\psi(k/2) e^{-iwk/2} \hat{S}^\phi(w/2) = Q_s(z)\hat{S}^\phi(w/2) \tag{33}$$

where

$$Q_s(z) = 1/2 \sum_{k\in\mathbb{Z}} S^\psi(k/2) e^{-iwk/2} = \sum_{k\in\mathbb{Z}} \hat{S}^\psi(w+4k\pi). \tag{34}$$

Since $\{V_j\}_{j\in\mathbb{Z}}$ forms an MRA and since $\{S^\phi(2^j x-k)\}_{k\in\mathbb{Z}}$ and $\{S^\psi(2^j x-k)\}_k$ respectively form Riesz bases for $V_j$ and $W_j$, $P_s(z)$ in (12) and $Q_s(z)$ in (34) should form a pair of reconstruction filters in $\{V_j\}_{j\in\mathbb{Z}}$ such that

$$\begin{cases} \hat{S}^\phi(w) = P_s(z)\hat{S}^\phi(w/2) \\ \hat{S}^\psi(w) = Q_s(z)\hat{S}^\phi(w/2) \end{cases}, \tag{35}$$

and

$$0 < A_\Delta \leq \left|\Delta_{P_s,Q_s}\right|^2 \leq B_\Delta < +\infty \tag{36}$$

(see Theorem 5.16 in [20]), where $A_\Delta$ and $B_\Delta$ are two constants and

$$\Delta_{P_s,Q_s} = P_s(z)Q_s(-z) - P_s(-z)Q_s(z). \tag{37}$$

Now, we discuss relationships between $P_s(z)$ and $Q_s(z)$. Furthermore, an algorithm is proposed for constructing $Q_s(z)$ directly from $S^\phi(x)$ and $P_s(z)$.

In order to explain our algorithm and the relationship between $Q_s(z)$ and $P_s(z)$, we first introduce

**Lemma 2.** Let $\psi(x)$ be the element of $W_0$ whose Fourier transform is expressed as

$$\hat{\psi}(w) = Q_\psi(z)\hat{\phi}(w/2), \tag{38}$$

with

$$Q_\psi(z) = \frac{1}{2}\sum_{k\in\mathbb{Z}} q_k z^k, \quad \{q_k\}_{k\in\mathbb{Z}} \in l^2(\mathbb{R}). \tag{39}$$

Then, we have

$$Q_\psi(z)E_\phi(z)\overline{P_\phi(z)} + Q_\psi(-z)E_\phi(-z)\overline{P_\phi(-z)} = 0, \tag{40}$$

where

$$E_\phi(z) = E_\phi(e^{-iw/2}) = \sum_{k=-\infty}^{+\infty} |\hat{\phi}(w/2 + 2k\pi)|^2, \quad z = e^{-iw/2}. \tag{41}$$

**Proof.** See Appendix E for the proof.

From Lemma 2, we have Theorem 2, which implies that the interpolation basis $\{S^\psi(2^j x - k)\}_k$ of $W_j$ has a close relationship with the interpolation basis $\{S^\phi(2^j x - k)\}_k$ of $V_j$ via the pair of reconstruction filters $(P_s(z), Q_s(z))$.

**Theorem 2** Assume that there exists an interpolation basis $\{S^\phi(x-k)\}_{k\in\mathbb{Z}}$ relative to the samples of $f_s(x)$ at $x = k \in \mathbb{Z}$ in the approximation space $V_0$. Then, if $W_j$ has an interpolation basis $\{S^\psi(2^j x - k)\}_{k\in\mathbb{Z}}$ relative to the samples of $f_s(x)$ at $x = 1/2^{j+1} + k/2^j$ with $k \in \mathbb{Z}$, it follows that

$$\begin{cases} P_s(-z) + P_s(z) = 1 \\ Q_s(z)/z - Q_s(-z)/z = 1 \end{cases} \tag{42}$$

and

$$Q_s(-1) = -1. \tag{43}$$

Furthermore, equation (23) holds and

$$Q_s(z) = \frac{zE_s(-z)\overline{P_s(-z)}}{P_s(-z)E_s(-z) + P_s(z)E_s(z)}. \tag{44}$$

**Proof.** See Appendix F for the proof.

When $S^\phi(x)$ is an interpolation scaling function and $S^\psi(x)$ is an interpolation wavelet corresponding to $S^\phi(x)$, $P_s(z)$ in (12) and $Q_s(z)$ in (34) form a pair of reconstruction filters. In this case, Theorem 2 indicates that (23) always holds. This implies that $A_s \leq |PE_s(w)|^2 \leq B_s$ is a necessary condition for the interpolation basis $\{S^\psi(2^j x-k)\}_{k\in\mathbb{Z}}$ to exist in $W_j$ when there exists the interpolation basis $\{S^\phi(x-k)\}_{k\in\mathbb{Z}}$. Hence, Theorems 1 and 2 imply that $A_s \leq |PE_s(w)|^2 \leq B_s$ is a necessary and sufficient condition for the existence of the interpolation basis $\{S^\psi(2^j x-k)\}_{k\in\mathbb{Z}}$ relative to $\{f_s(k/2^j+1/2^{j+1})\}_{k\in\mathbb{Z}}$ in $W_j$ when there exists the interpolation basis $\{S^\phi(x-k)\}_{k\in\mathbb{Z}}$.

Furthermore, Theorem 2 shows that $Q_s(z)$ can be connected to $P_s(z)$ by (44). This implies that $S^\phi(x)$ and $S^\psi(x)$ have a close relationship in their analytic expressions. Hence, we can construct the interpolation wavelet $S^\psi(2^j x)$ in $W_j$ from the interpolation scaling function $S^\phi(x)$.

In addition, since the interpolation function $S^\phi(x)$ can be constructed from a general scaling function $\phi(x)$ as in (10), $E_s(z)$ and $P_s(z)$ in (44) can be constructed from a scaling function $\phi(x)$. Hence, (33) and (44) imply that we can construct the interpolation function $S^\psi(x)$ using any scaling function $\phi(x)$ when $S^\psi(x)$ and $S^\phi(x)$ exist.

Finally, for a pair of reconstruction filters $\left(P_\phi(z), Q_\psi(z)\right)$, we have the equality

$$\begin{cases} P_\phi(1)=1 \quad P_\phi(-1)=0 \\ Q_\psi(1)=0 \end{cases} \tag{45}$$

(see equation (5.4.1) in [26]), which is used widely in numerical computations of wavelet theory. Obviously, since $\{S^\phi(x-k)\}_{k\in\mathbb{Z}}$ and $\{S^\psi(x-k)\}_{k\in\mathbb{Z}}$ are respectively Riesz bases of $V_0$ and $W_0$, equation (45) also holds for $\left(P_s(z), Q_s(z)\right)$. From (45), $P_\phi(z)$ has values at $z=\pm 1$ that are independent of a scaling function $\phi(x)$, but this holds for $Q_\psi(z)$ only at $z=1$. However, by (43) and (45), $Q_s(z)$ again has such independent values at $z=\pm 1$. Hence, besides (45), equation (43) is also an useful condition for numerical computations of an interpolation function $S^\psi(x)$ in a multiresolution analysis.

**Remark.** Theorem 2 indicates that $A_s \leq |PE_s(w)|^2 \leq B_s$ is a necessary condition for the interpolation basis $\{S^\psi(2^j x-k)\}_{k\in\mathbb{Z}}$ to exist in $W_j$, when there exists the interpolation basis $\{S^\phi(x-k)\}_{k\in\mathbb{Z}}$. Simultaneously, it shows that $Q_s(z)$ has the expression in terms of $P_s(z)$ and $E_s(z)$ as in (44). Since $P_s(z)$ and $E_s(z)$ can be constructed from $S^\phi(x)$ by (12) and (16), and since $S^\phi(x)$ can be constructed from a scaling function $\phi(x)$ by (10), equations (33) and (44) imply that the interpolation function $S^\psi(x)$ in a wavelet space can be constructed from a scaling function $\phi(x)$. Hence, Theorem 2 has supplied an effective method for constructing an interpolation wavelet $S^\psi(x)$ in $W_0$. Finally, since $\hat{S}^\psi(w)$ has an analytic expression involving $\hat{S}^\phi(w)$ as in (33) and (44), the interpolation basis $\{S^\psi(x-k)\}_{k\in\mathbb{Z}}$ in a wavelet space is closely related to the interpolation basis $\{S^\phi(x-k)\}_{k\in\mathbb{Z}}$ in an approximation space.

## 5. SIMULATION

### 5.1 Wavelet Spaces in Simulation

Cardinal splines of different orders can be used to form many important interpolation functions, for example the Haar or triangular scaling functions. Note also that interpolation bases for cardinal spline spaces have varying smoothness and orthogonality properties, depending on spline orders. It follows that cardinal spline wavelets can form a good experimental testbed for the theories of this paper, as is shown below. Additionally, since the Shannon wavelet is a typical interpolation wavelet, we also include it in our simulation below.

In the simulation, the wavelet spaces are divided into two groups. The Shannon wavelet and the Haar wavelet spaces form the first group, for which the interpolation wavelets are known from the standard wavelet theory. The spaces spanned by the Cardinal splines of orders 2, 3 and 4 form the second group, for which the interpolation wavelets are not widely known.

For the first group, we will construct the interpolation wavelets using an algorithm based on Theorem 2. Then, these interpolation wavelets constructed by our algorithm are compared to their counterparts made from standard wavelet constructions.

For the second group, we will first obtain a standard wavelet $\psi(x)$ whose Fourier transform is

$$\hat{\psi}(w) = -z E_\phi(-z) \overline{P_\phi(-z)} \hat{\phi}(w/2) \qquad (46)$$

(equation (5.6.13) in [20]), where $\hat{\phi}(w)$ denotes the Fourier transform of the scaling function $\phi(x)$ in the approximation space $V_0$, and $P_\phi(z)$ and $E_\phi(z)$ are respectively defined in (9) and (41). Then the interpolation functions $\{S^\psi(x-k)\}_{k\in\mathbb{Z}}$ constructed by our algorithm will be applied to interpolate the above standard wavelet $\psi(x)$, as indicated below. By equation (5.6.13) in [20], the wavelets $\{\psi(x-k)\}_{k\in\mathbb{Z}}$ in (46) form a Riesz basis for the wavelet space $W_0$. Hence, if our interpolation functions $\{S^\psi(x-k)\}_{k\in\mathbb{Z}}$ can recover the wavelet $\psi(x)$ in (46) from its samples, then $\{S^\psi(x-k)\}_{k\in\mathbb{Z}}$ also form a Riesz basis in the corresponding wavelet space $W_0$. This verifies that $S^\psi(x)$ is an interpolation wavelet for $W_0$.

## *5.2 Shannon and Haar Wavelets*

The Shannon and the Haar scaling functions are two well-known interpolation functions that form orthogonal bases in their respective spaces. In this section, we obtain the interpolation wavelets corresponding to them by using (33) and (44). More specifically, the interpolation wavelets obtained from our algorithm will be compared to the Shannon and the Haar wavelets, so that we can verify our theorems for the Shannon and the Haar MRA.

### *5.2.1 Existence of interpolation basis in the Shannon and the Haar wavelet space*

The scaling function $S^\phi(x) = \sin(\pi x)/(\pi x)$ for the Shannon MRA is depicted in Figure 1-a, and the scaling function $S^\phi(x) = 1_{[0,1]}(x)$ for the Haar MRA is in Figure 1-b. Note that either the Shannon or the Haar scaling functions can form an interpolation basis $\{S^\phi(x-k)\}_{k\in\mathbb{Z}}$ relative to the samples $\{f_s(k)\}_{k\in\mathbb{Z}}$ in the approximation spaces $V_0$ for their respective multiresolution analysis. On the other hand, both the Shannon and the Haar wavelet are known as interpolation functions relative to $\{f_s(k+1/2)\}_{k\in\mathbb{Z}}$ in $W_0$ [6].

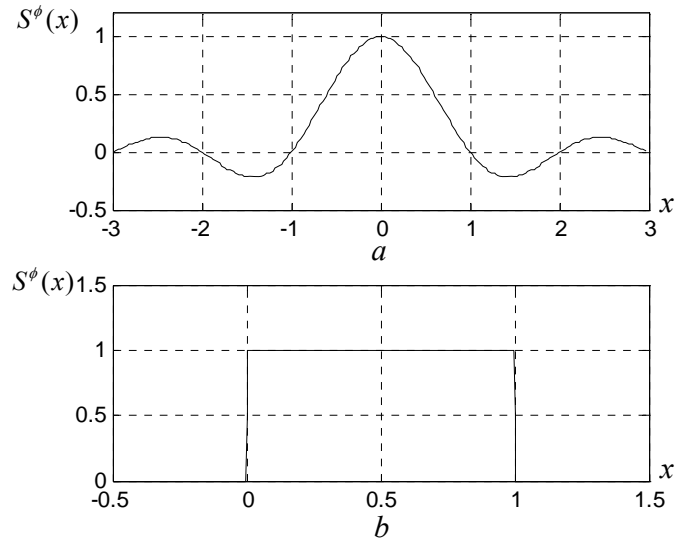

Figure 1. Scaling function $S^\phi(x)$. (a) The Shannon scaling function. (b) The Haar scaling function.

Hence, in this section, we verify that (23) is a necessary condition for an interpolation wavelet to exist by showing that (23) holds for the Shannon and the Haar scaling functions.

It follows from (12) that

$$P_s(z) = \sum_{k=-\infty}^{+\infty} 1_{[-\pi,\pi]}(w+4k\pi) \quad (47)$$

for the Shannon scaling function and

$$P_s(z) = 1/2 \times \sin(w/2)e^{-iw/2} \times \lim_{\varepsilon \to 0^+} \sum_{k=-\infty}^{+\infty} e^{-i(w+4k\pi)\varepsilon}/(w/4+k\pi)$$

$$= 1/2 \times \sin(w/2)e^{-iw/2} \times (\sum_{k=-\infty}^{+\infty} 1/(w/4+k\pi) + 4i\int_{-\infty}^{+\infty} \frac{\sin(4\pi w)}{4\pi w} dw)$$

$$= 1/2 \times \sin(w/2)e^{-iw/2}(\cot(w/4)+i) = (1+e^{-iw/2})/2 \quad . \quad (48)$$

for the Haar scaling function.

Since the Shannon and Haar scaling functions, both denoted as $\{S^\phi(x-k)\}_{k\in\mathbb{Z}}$ here, form the orthogonal bases for $V_0$, we have

$$E_s(z) = 1 \quad (49)$$

for both scaling functions.

Inserting (47) and (49) into (15) yields

$$PE_s(w) = \overline{\sum_{k=-\infty}^{+\infty} 1_{[-\pi,\pi]}(w+2\pi+4k\pi)} + \overline{\sum_{k=-\infty}^{+\infty} 1_{[-\pi,\pi]}(w+4k\pi)}$$

$$= \sum_{k=-\infty}^{+\infty} 1_{[-\pi,\pi]}(w+2k\pi) = 1 \quad (50)$$

for the Shannon scaling function and inserting (48) and (49) into (15) yields

$$PE_s(w) = \overline{(1+e^{-iw/2})/2} + \overline{(1+e^{-i(w+2\pi)/2})/2} = 1/2 + e^{iw/2}/2 + 1/2 - e^{iw/2}/2 = 1. \quad (51)$$

Equations (50) and (51) indicate that $|PE_s(w)|^2 = 1$ for both the Shannon and the Haar scaling functions, i.e. equation (23) holds. This verifies that (23) is a necessary condition for the existence of interpolation wavelets in the Shannon and the Haar wavelet space.

### 5.2.2  *Constructing interpolation wavelets*

From Theorem 2, an interpolation wavelet $S^\psi(x)$ can be constructed from an interpolation scaling function $S^\phi(x)$. Here, we construct the interpolation wavelets $S^\psi(x)$ corresponding to the Shannon and the Haar scaling function respectively by using (33) and (44). Then, we verify our constructions by comparing them to their respective (the Shannon and the Haar) wavelets.

Since each of these two scaling functions forms an interpolation basis $\{S^\phi(x-k)\}_{k\in\mathbb{Z}}$ relative to the samples $\{f_s(k)\}_{k\in\mathbb{Z}}$, the interpolation functions $S^\psi(x)$ in the corresponding wavelet spaces $W_0$ can be obtained directly from (33) and (44) based on Theorem 2.

Inserting (50) or (51) into (44) yields

$$Q_s(z) = zE_s(-z)\overline{P_s(-z)}. \qquad (52)$$

for both the Shannon and the Haar MRAs. Since (49) indicates $E_s(z) = 1$ for both the Shannon and the Haar scaling function, equation (52) is rewritten as

$$Q_s(z) = z\overline{P_s(-z)}. \qquad (53)$$

Inserting (47) and (48) into (53) yields respectively

$$Q_s(z) = z \times \overline{\sum_{k=-\infty}^{+\infty} 1_{[-\pi,\pi]}(w + 2\pi + 4k\pi)} \qquad (54)$$

for the Shannon MRA and

$$Q_s(z) = \frac{1}{2}e^{-iw/2} - \frac{1}{2} \qquad (55)$$

for the Haar MRA.

Inserting (54) and (55) in (33) and taking the inverse Fourier transform yields the interpolation function

$$S^\psi(x) = \frac{\sin\frac{\pi}{2}(x - 1/2)}{\frac{\pi}{2}(x - 1/2)} \cos\frac{3\pi}{2}(x - 1/2) \qquad (56)$$

for the Shannon wavelet space, and the interpolation function

$$S^\psi(x) = -1_{[0,\frac{1}{2}]}(x) + 1_{[\frac{1}{2},1]}(x). \qquad (57)$$

for the Haar wavelet space. The interpolation wavelets in (56) and (57) are respectively described in Figures 2-a and b.

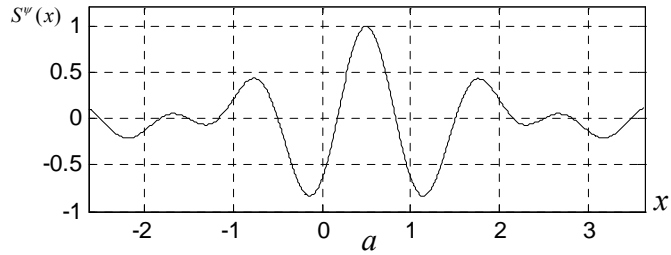

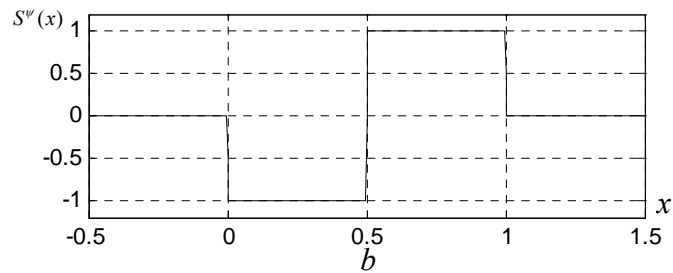

Figure 2. Interpolation functions $S^\psi(x)$. (a) $S^\psi(x)$ in Shannon wavelet space. (b) $S^\psi(x)$ in Haar wavelet space.

By comparing the interpolation wavelets $S^\psi(x)$ in (56) and (57) to the Shannon and Haar wavelets, it follows that the interpolation wavelets $S^\psi(x)$ in (56) and (57) are the same as their counterparts from the standard theory. Hence, the simulations of this section have validated our algorithm for the Shannon and the Haar wavelet spaces.

## 5.3 The spline function of order 2

### 5.3.1. Existence of interpolation function in $W_0$

The spline scaling function of order 2 (known as the triangular function) can be considered a typical interpolation function, and is depicted in Figure 3-a. Its corresponding Fourier transform is

$$\hat{S}^\phi(w) = \left(\frac{\sin(w/4)}{w/4}\right)^2. \tag{58}$$

Since the translates of the triangular functions form an interpolation basis $\{S^\phi(x-k)\}_k$ relative to $\{f_s(k)\}_{k\in\mathbb{Z}}$ in $V_0$, Theorems 1 and 2 imply that (23) is a necessary and sufficient condition for an interpolation wavelet $S^\psi(x)$ to exist. Hence, in this section, we first use (23) to determine existence of an interpolation spline wavelet of order 2.

It follows from (12) and (58) that the two-scale symbol of the spline scaling function of order 2 is

$$P_s(z) = e^{-iw/2}/4 + 1/2 + e^{iw/2}/4. \tag{59}$$

Inserting (59) into (15) yields

$$PE_s(w) = \overline{P_s(-z)}E_s(-z) + \overline{P_s(z)}E_s(z) = (2 + \cos^2(w/2))/3. \tag{60}$$

Obviously, equation (60) indicates $4/9 \leq |PE_s(w)|^2 \leq 1$ for the spline scaling function of order 2, which implies that (23) holds for the spline function of order 2. Hence, there exists an interpolation basis $\{S^\psi(x-k)\}_k$ relative to $\{f_s(k+1/2)\}_{k\in\mathbb{Z}}$ in $W_0$ according to Theorem 1 and Proposition 3.

### 5.3.2. Construction of interpolation function in $W_0$

By Theorem 1 and Proposition 3, equation (60) indicates that there exists an interpolation spline wavelet $S^\psi(x)$ of order 2, so we construct $S^\psi(x)$ using (33) and (44) in this section.

Inserting (59) and (60) into (44) yields

$$Q_s(z) = \frac{e^{-iw/2}\sin^2(w/4)(1+2\sin^2(w/4))}{2+\cos^2(w/2)}. \tag{61}$$

It follows from (33) and (61) that the Fourier transform of the interpolation wavelet $S^\psi(x)$ is

$$\hat{S}^\psi(w) = \frac{e^{-iw/2}\sin^2(w/4)(1+2\sin^2(w/4))}{2+\cos^2(w/2)} \times \left(\frac{\sin(w/4)}{w/4}\right)^2. \tag{62}$$

The corresponding interpolation function $S^\psi(x)$ is depicted in Figure 3-b, which is obtained by calculating numerically the inverse Fourier transform of $\hat{S}^\psi(w)$ in (62).

### 5.3.3. Illustration of wavelet recovery

Since $\psi(x)$ in (46) is a wavelet, we verify that $S^\psi(x)$ of Figure 3-b is an interpolation wavelet by using $\{S^\psi(x-k)\}_{k\in\mathbb{Z}}$ to interpolate $\psi(x)$.

It follows from (46) that the Fourier transform of a spline wavelet $\psi(x)$ of order 2 is

$$\hat{\psi}(w) = \frac{e^{-iw/2} - 6 + 10e^{iw/2} - 6e^{iw} + e^{3iw/2}}{24}\left(\frac{\sin(w/4)}{w/4}\right)^2. \tag{63}$$

The spline wavelet $\psi(x)$ shown in Figure 4-a is obtained by calculating numerically the inverse Fourier transform of $\hat{\psi}(w)$ in (63). The interpolation functions $\{S^\psi(x-k)\}_{k\in\mathbb{Z}}$ with $S^\psi(x)$ depicted in Fig. 3-b are applied to interpolate the samples $\{\psi(n-1/2)\}_{n\in\mathbb{Z}}$, $n=-3,\cdots,3$. The recovery $f_{ap}(x)$ of $\psi(x)$ is shown in Fig. 4-b. The error $e(x)$ of $f_{ap}(x)$ is depicted in Figure 4-c. Since both $S^\psi(x)$ and $\psi(x)$ are obtained from the numerical calculation in our simulation, the accuracy of recovery is connected with the quantity of data in the numerical calculation. Hence, an error $e(x)$ between the recovery $f_{ap}(x)$ and the wavelet $\psi(x)$ cannot be avoided because it arises from computer memory limitations. However, it is known from Figure 4-c that the error $e(x)$ of recovery $f_{ap}(x)$ is smaller than $2\times 10^{-4}$, which is reasonable for our numerical calculation.

Hence, we can conclude that $\psi(x)$ has been recovered successfully using our interpolation functions $\{S^\psi(x-k)\}_{k\in\mathbb{Z}}$. This result indicates that the constructed interpolation functions $\{S^\psi(x-k)\}_{k\in\mathbb{Z}}$ can form an interpolation basis in $W_0$, and so $S^\psi(x)$ is an interpolation wavelet. This experiment verifies Theorem 1 and Proposition 3 for the spline function of order 2, i.e., equation (23) is a sufficient condition for the existence of the interpolation spline wavelet of order 2. On the other hand, this experiment also shows that our algorithm in Theorem 2 can construct an interpolation wavelet from an interpolation scaling function.

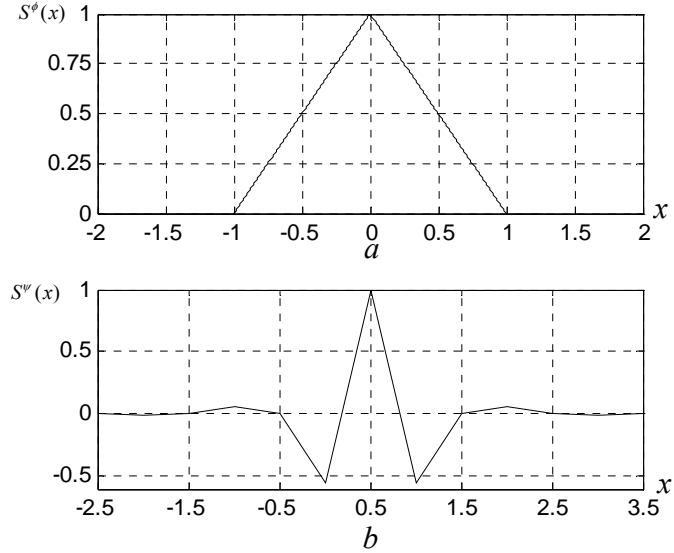

Figure 3. Spline scaling function $S^\phi(x)$ of order 2 and its corresponding interpolation function $S^\psi(x)$ in $W_0$. (a) Scaling function $S^\phi(x)$. (b) Interpolation function $S^\psi(x)$ in $W_0$.

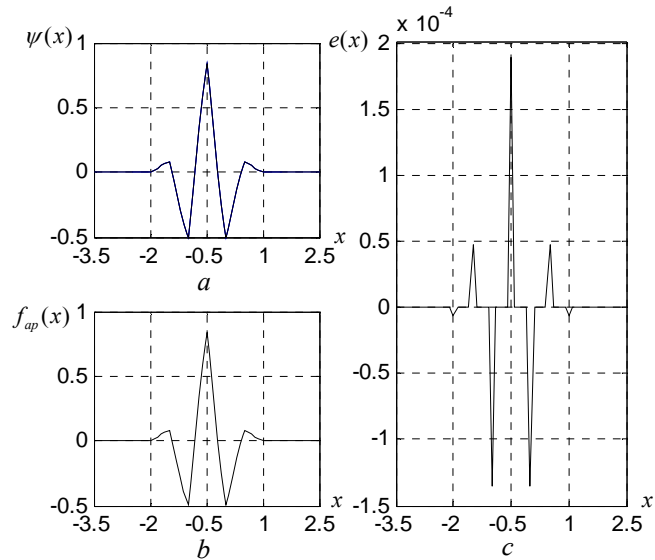

Figure 4. Spline wavelet $\psi(x)$ of order 2 and its recovery $f_{ap}(x)$. (a) Spline wavelet $\psi(x)$ of order 2. (b) Its recovery $f_{ap}(x)$. (c) Error $e(x)$ of recovery.

## 5.4 *The spline function of order 3*

The spline function $\phi(x)$ of order 3 is shown in Figure 5 and its Fourier transform is

$$\hat{\phi}(w) = \left(\frac{\sin(w/2)}{(w/2)}\right)^3. \tag{64}$$

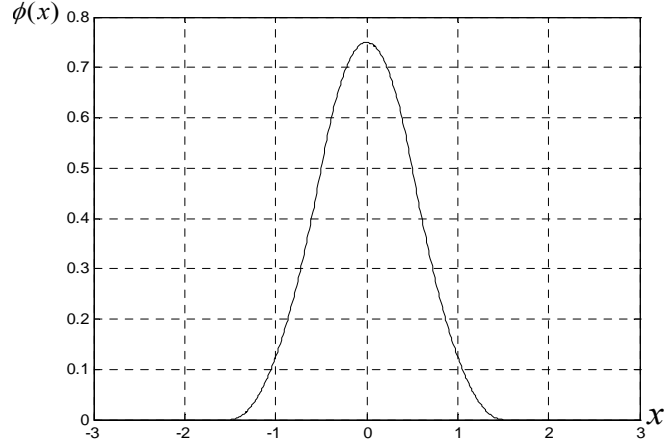

Figure 5. Spline function $\phi(x)$ of order 3.

For this spline, we have

$$1/2 \leq \sum_{k=-\infty}^{+\infty} \hat{\phi}(w+2k\pi) = \frac{\cos^2(w/2)+1}{2} \leq 1. \tag{65}$$

By Theorem 1 in [8], an interpolation basis $\{S^\phi(x-k)\}_{k\in\mathbb{Z}}$ exists in an approximation space if and only if $1/\sum_{k=-\infty}^{+\infty}\hat{\phi}(w+2k\pi) \in L^2[0,2\pi]$. Hence, equation (65) indicates that there exists a basis $\{S^\phi(x-k)\}_{k\in\mathbb{Z}}$ in the third order spline space.

Inserting (64) into (10), we obtain the Fourier transform of $S^\phi(x)$ as

$$\hat{S}^\phi(w) = \frac{16\sin^3(w/2)}{w^3(\cos^2(w/2)+1)}. \tag{66}$$

By taking the inverse Fourier transform of (66), we obtain the interpolation function $S^\phi(x)$, shown in Figure 6-a, for the third order spline space.

Since there exists an interpolation basis $\{S^\phi(x-k)\}_{k\in\mathbb{Z}}$ in the third order spline space, there exists an interpolation basis $\{S^\psi(x-k)\}_{k\in\mathbb{Z}}$ relative to $\{f_s(k+1/2)\}_{k\in\mathbb{Z}}$ in the corresponding wavelet space $W_0$ if and only if (23) holds, by Theorems 1 and 2. Hence, we apply (23) to decide existence of an interpolation wavelet $S^\psi(x)$.

Inserting (66) into (15) yields

$$PE_s(w) = \overline{P_s(-z)}E_s(-z) + \overline{P_s(z)}E_s(z) = \frac{4}{15} \times \frac{\cos(w/4) - \sin(w/4)}{1+\cos^2(w/2)}$$
$$\times \left( \frac{15 + 15\sin(w/4)\cos(w/4)}{(1+\cos^2(w/4))(1+\sin^2(w/4))} + \frac{11\sin^3(w/4)\cos^3(w/4) - 2\sin^4(w/4)\cos^4(w/4)}{(1+\cos^2(w/4))(1+\sin^2(w/4))} \right), \quad (67)$$

which is depicted in Figure 6-b. This figure and equation (67) indicate that $PE_s(w) = 0$ for $w = 4k\pi + \pi$, which implies that (23) does not hold. By Theorem 2, if there exists an interpolation basis $\{S^\psi(x-k)\}_{k\in\mathbb{Z}}$ in $W_0$, then (23) must hold for all $w \in \mathbb{R}$. Hence, the condition $PE_s(w) = 0$ for $w = 4k\pi + \pi$ in (67) implies that there does not exist an interpolation basis $\{S^\psi(x-k)\}_{k\in\mathbb{Z}}$ in the wavelet space $W_0$.

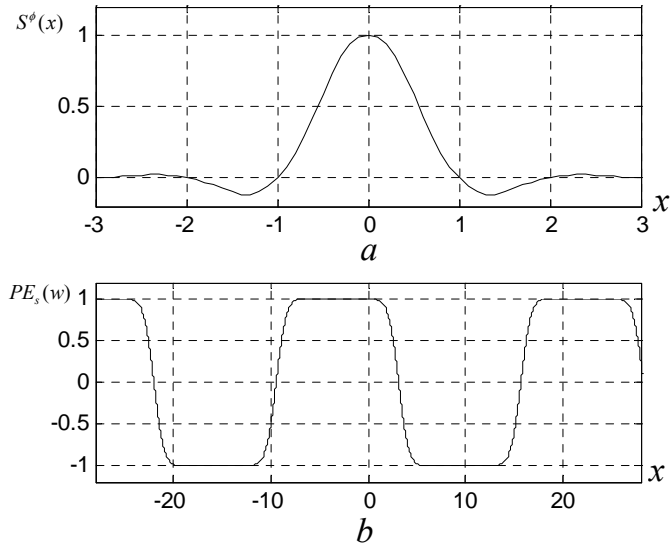

Figure 6. Interpolation function $S^\phi(x)$ of the third order spline space and its corresponding $PE_s(w)$. (a) Interpolation scaling function $S^\phi(x)$. (b) Function $PE_s(w)$.

In fact, the spline function $\phi(x)$ with Fourier transform (64) cannot be represented as a linear combination of $\{\phi(2x-k)\}_k$, so the spline function of order 3 is not a scaling function in this case. Hence, there does not exist a wavelet space corresponding to the third order spline space, so there also does not exist an interpolation wavelet $S^\psi(x)$ corresponding to $S^\phi(x)$.

Hence, the above numerical simulation for the third order spline verifies that (23) is a necessary condition for the existence of an interpolation wavelet.

## 5.5 The spline function of order 4

### 5.5.1. Existence of interpolation function in $W_0$

Figure 7-a depicts the spline of order 4 whose Fourier transform is represented as

$$\hat{\phi}(w) = \left(\sin(w/2)/(w/2)\right)^4. \tag{68}$$

Obviously, this spline is a scaling function (equation (7.11) in [6]). From (68), we have

$$1/3 \le \sum_{k=-\infty}^{+\infty} \hat{\phi}(w+2k\pi) = \frac{1}{3}(2\cos^2(w/2)+1) \le 1. \tag{69}$$

Equation (69) indicates that there exists an interpolation basis $\{S^\phi(x-k)\}_{k\in\mathbb{Z}}$ in the space $V_0$ spanned by splines of order 4 (see Theorem 1 in [8]). Figure 7-b shows the spline interpolation function $S^\phi(x)$ of order 4 with Fourier transform

$$\hat{S}^\phi(w) = \frac{3}{2\cos^2(w/2)+1}\left(\frac{\sin(w/2)}{w/2}\right)^4. \tag{70}$$

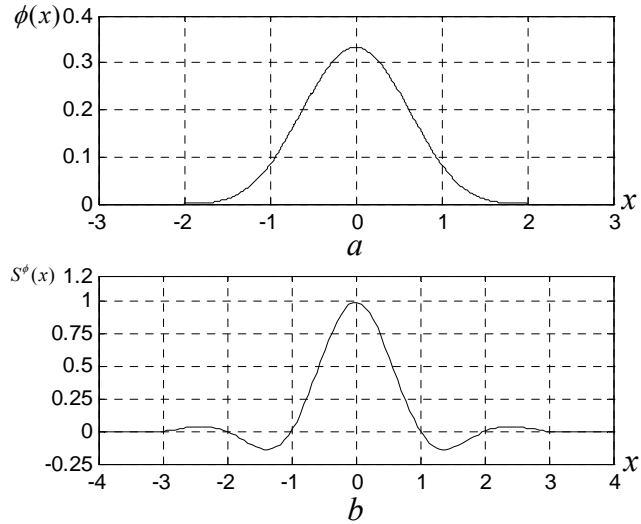

Figure 7. Spline scaling function $\phi(x)$ of order 4 and its corresponding interpolation scaling function $S^\phi(x)$. (a) Scaling function $\phi(x)$. (b) Interpolation scaling function $S^\phi(x)$.

Since there exists an interpolation basis $\{S^\phi(x-k)\}_k$ relative to $\{f_s(k)\}_{k\in\mathbb{Z}}$ in the fourth order spline space, equation (23) is a necessary and sufficient condition for an interpolation wavelet $S^\psi(x)$ to exist in the corresponding wavelet space $W_0$, by Theorems 1 and 2. Hence, we apply (23) to determine existence of interpolation spline wavelet $S^\psi(x)$ of order 4.

It follows from (12) and (70) that

$$P_s(z) = \sum_{k \in \mathbb{Z}} \hat{S}^\phi(w + 4k\pi) = \frac{\cos^4(w/4)}{2\cos^2(w/2) + 1}\left(2\cos^2(w/4) + 1\right). \quad (71)$$

From (15), (70) and (71), we have

$$PE_s(w) = \overline{P_s(-z)}E_s(-z) + \overline{P_s(z)}E_s(z) = \frac{3\cos(w)}{2240} - \frac{129}{2240(\cos(w)+2)} + \frac{39}{140(\cos(w)-7)} + 149/140, \quad (72)$$

which is depicted in Figure 8-a.

Figure 8-a illustrates the fact that $0.9705 < |PE_s(w)|^2 < 1$. From Theorem 1 and Proposition 3, it follows that there exists an interpolation function $\{S^\psi(x-k)\}_k$ relative to $\{f_s(k+1/2)\}_{k \in \mathbb{Z}}$ in $W_0$.

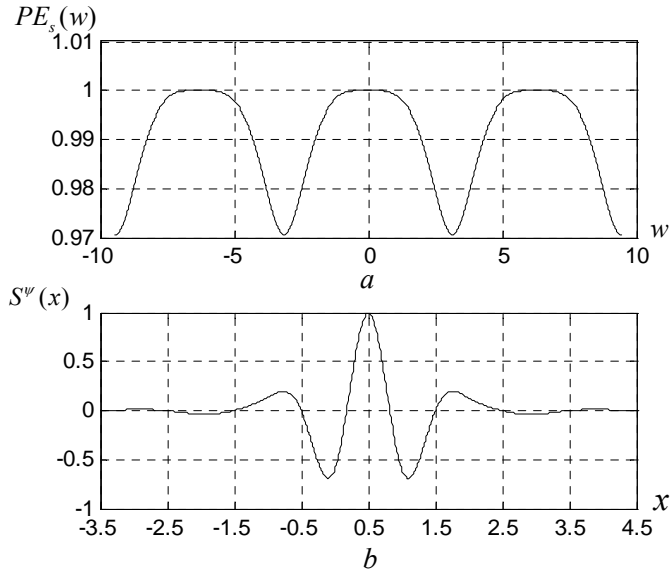

Figure 8. $PE_s(w)$ and interpolation wavelet $S^\psi(x)$ corresponding to the spline function of order 4 (a) Function $PE_s(w)$. (b) Interpolation wavelet $S^\psi(x)$.

### 5.5.2. *Construction of interpolation wavelet*

By Theorem 1 and Proposition 3, equation (72) indicates that there exists an interpolation spline wavelet $S^\psi(x)$ of order 4, so we construct $S^\psi(x)$ by using (33) and (44) in this section.

Since equation (23) holds, inserting (70) into (44) yields

$$Q_s(z) = \frac{e^{-iw/2} \times (\cos(w/2) - 1)^2}{-2\sin^6(w/2) + 594\sin^4(w/4) + 840\sin^2(w/2) - 2520} \\ \times \left(\cos^4(w/4) - 58\cos^3(w/2) + 177\cos^2(w/2) + 322\cos(w/2) - 544\right). \quad (73)$$

Then, by using (33) and (70), we have the Fourier transform of $S^\psi(x)$ as

$$\hat{S}^\psi(w) = Q_s(z) \times \frac{3}{2\cos^2(w/4) + 1}\left(\frac{\sin(w/4)}{w/4}\right)^4 \quad (74)$$

with $Q_s(z)$ in (73).

By calculating numerically the inverse Fourier transform of $\hat{S}^\psi(w)$ in (74), we obtain the interpolation $S^\psi(x)$, depicted in Figure 8-b.

### 5.5.3. Illustration of wavelet recovery

Since $\psi(x)$ in (46) is a wavelet, we verify that $S^\psi(x)$ of Figure 8-b is an interpolation wavelet by using $\{S^\psi(x-k)\}_{k \in \mathbb{Z}}$ to interpolate $\psi(x)$.

It follows from (46) that

$$\hat{\psi}(w) = \left( \frac{1}{80640} e^{-3iw} - \frac{31}{20160} e^{-5iw/2} + \frac{559}{26880} e^{-2iw} - \frac{247}{2520} e^{-3iw/2} + \frac{9241}{40320} e^{-iw} - \frac{337}{1120} e^{-iw/2} \right.$$
$$\left. + \frac{9241}{40320} - \frac{247}{2520} e^{iw/2} + \frac{559}{26880} e^{iw} - \frac{31}{20160} e^{3iw/2} + \frac{1}{80640} e^{2iw} \right) \times \left( \frac{\sin(w/4)}{w/4} \right)^4 \quad (75)$$

is the Fourier transform of a spline wavelet $\psi(x)$ of order 4. Figure 9-a shows the wavelet $\psi(x)$ obtained by calculating numerically the inverse Fourier transform of $\hat{\psi}(w)$ in (75). The samples $\{\psi(n-1/2)\}_{n \in \mathbb{Z}}$ at $n = -4, \cdots, 4$, are used to recover $\psi(x)$. The recovery $f_{ap}(x)$ of $\psi(x)$ is shown in Figure 9-b. Since $S^\psi(x)$ and $\psi(x)$ in our simulation are obtained from numerical calculations, the recovery cannot avoid an error $e(x)$, shown in Figure 9-c. From this figure, it is seen that the error $e(x)$ is smaller than $10 \times 10^{-5}$. This error level is reasonable for the recovery in our algorithm.

Hence, this simulation shows that $\psi(x)$ has been recovered, verifying that $\{S^\psi(x-k)\}_{k \in \mathbb{Z}}$ obtained by our algorithm is the interpolation basis relative to $\{f_s(k+1/2)\}_{k \in \mathbb{Z}}$. Hence, we conclude that equation (23) is a sufficient condition for the existence of the interpolation spline wavelet of order 4.

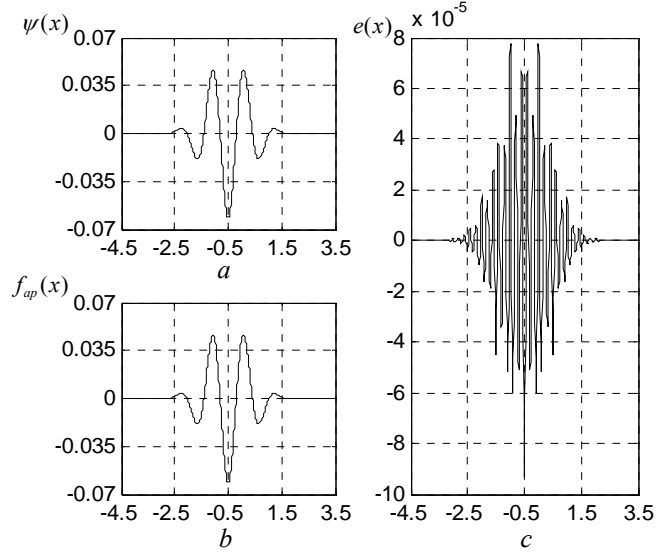

Figure 9. Spline wavelet $\psi(x)$ of order 4 and its recovery $f_{ap}(x)$. (a) Spline wavelet $\psi(x)$ of order 4. (b) Its recovery $f_{ap}(x)$. (c) Error $e(x)$ of recovery.

## 5.6 Simulation Summary

The simulations in Section 5.2 illustrate that equation (23) ( $0 < A_s \leq |PE_s(w)|^2 \leq B_s < +\infty$ ) holds for the Shannon and Haar scaling functions due to the fact their corresponding wavelets are interpolation functions. Furthermore, it is demonstrated in Section 5.4 that (23) fails to hold when there does not exist a wavelet space corresponding to the third order spline space. On the other hand, the simulations in Sections 5.3 and 5.5 illustrate that we can always construct an interpolation wavelet $S^\psi(x)$ when (23) holds for the spline functions of orders 2 and 4. Hence, our simulation results validate experimentally the fact that (23) is a necessary and sufficient condition for the existence of an interpolation basis $\{S^\psi(2^j x - k)\}_{k \in \mathbb{Z}}$ relative to $\{f_s(1/2^{j+1} + k/2^j)\}_{k \in \mathbb{Z}}$ in $W_j$.

In addition, the simulations in Section 5.2 illustrate that the algorithm based on Theorem 2 has succeeded in reconstructing the Shannon and Haar wavelets. Simultaneously, in Sections 5.3 and 5.5, the algorithm based on Theorem 2 is applied to construct the interpolation wavelets $\{S^\psi(x-k)\}_{k \in \mathbb{Z}}$ corresponding to the spline functions of orders 2 and 4. By interpolating the wavelet $\psi(x)$ in (46), we have verified that these constructed functions $\{S^\psi(x-k)\}_{k \in \mathbb{Z}}$ are the interpolation spline wavelets of orders 2 and 4. Hence, the simulations in this paper also validate that the algorithm based on Theorem 2 can indeed construct an interpolation wavelet from an interpolation scaling function.

## 6. CONCLUSION

In an MRA, interpolation bases in wavelet spaces $\{W_j\}_{j\in\mathbb{Z}}$ and approximation spaces $\{V_j\}_{j\in\mathbb{Z}}$ are related closely to each other. We have verified that if there exists an interpolation basis $\{S^\phi(x-k)\}_{k\in\mathbb{Z}}$ relative to $\{f_s(k)\}_{k\in\mathbb{Z}}$ in $V_0$, then $0 < A_s \leq \left|\overline{P_s(-z)}E_s(-z) + \overline{P_s(z)}E_s(z)\right|^2 \leq B_s < +\infty$ is a necessary and sufficient condition for the existence of the interpolation basis $\{S^\psi(2^j x - k)\}_{k\in\mathbb{Z}}$ relative to $\{f_s(k/2^j + 1/2^{j+1})\}_{k\in\mathbb{Z}}$ in all wavelet spaces $W_j$. In addition, it has been shown that $S^\phi(x)$ is a scaling function (an interpolation scaling function) for $V_0$, so $(P_s(z), Q_s(z))$ form a pair of reconstruction filters for $\{V_j\}_{j\in\mathbb{Z}}$. It has been indicated in Theorem 2 that $(P_s(z), Q_s(z))$ have a number of interesting properties and are connected to each other in the form of an analytic expression. With this connection, we propose an algorithm for constructing an interpolation wavelet $S^\psi(x)$ from an interpolation scaling function $S^\phi(x)$ directly. These results above all indicate that an interpolation basis $\{S^\psi(2^j x - k)\}_{k\in\mathbb{Z}}$ in $W_j$ is related closely to an interpolation basis $\{S^\phi(2^j x - k)\}_{k\in\mathbb{Z}}$ in $V_j$. The simulation indicates that our theorems are effective for determining the existence of $S^\psi(x)$, and for constructing $S^\psi(x)$ from an interpolation scaling function $S^\phi(x)$.

**Appendix A PROOF OF PROPOSITION 1**

**Part 1: Proof that the series $\sum_{k\in\mathbb{Z}} S^\gamma(k/2)e^{-iwk/2}$ and $\sum_{k\in\mathbb{Z}} S^\gamma(w+4k\pi)$ converge pointwise everywhere, so that $\frac{1}{2}\sum_{k\in\mathbb{Z}} S^\gamma(k/2)e^{-iwk/2} = \sum_{k\in\mathbb{Z}} S^\gamma(w+4k\pi)$ holds by the Poisson summation formula.**

Since the series $\sum_{k\in\mathbb{Z}} \hat{\gamma}(w+4k\pi)$ and $\sum_{k\in\mathbb{Z}} \gamma(k/2)e^{-iwk/2}$ both converge, we have

$$\sum_{k\in\mathbb{Z}} \hat{\gamma}(w+4k\pi) = \frac{1}{2}\sum_{k\in\mathbb{Z}} \gamma(k/2)e^{-iwk/2} \tag{A1}$$

by the Poisson summation formula.

In fact, since

$$\begin{cases} \sum_{k\in\mathbb{Z}} \hat{\gamma}(w+2k\pi) = \sum_{k\in\mathbb{Z}} \hat{\gamma}(w+4k\pi+2\pi) + \sum_{k\in\mathbb{Z}} \hat{\gamma}(w+4k\pi) \\ \sum_{k\in\mathbb{Z}} \hat{\gamma}(k)e^{-iwk} = \frac{1}{2}\sum_{k\in\mathbb{Z}} \hat{\gamma}(k/2)e^{-i(w+2\pi)k/2} + \frac{1}{2}\sum_{k\in\mathbb{Z}} \hat{\gamma}(k/2)e^{-iwk/2} \end{cases} \tag{A2}$$

(or

$$\begin{cases} \sum_{k\in\mathbb{Z}} (-1)^k \hat{\gamma}(w+2k\pi) = \sum_{k\in\mathbb{Z}} \hat{\gamma}(w+4k\pi) - \sum_{k\in\mathbb{Z}} \hat{\gamma}(w+4k\pi+2\pi) \\ e^{-iw/2}\sum_{k\in\mathbb{Z}} \hat{\gamma}(k+1/2)e^{-iwk} = \frac{1}{2}\sum_{k\in\mathbb{Z}} \hat{\gamma}(k/2)e^{-iwk/2} - \frac{1}{2}\sum_{k\in\mathbb{Z}} \hat{\gamma}(k/2)e^{-i(w+2\pi k)/2} \end{cases}, \tag{A3}$$

it is derived from (A2) (or from (A3)) that the series $\sum_{k\in\mathbb{Z}} \gamma(k)e^{-iwk}$ (or $\sum_{k\in\mathbb{Z}} \gamma(k+1/2)e^{-iwk}$) and

$\sum_{k\in\mathbb{Z}} \hat{\gamma}(w+2k\pi)$ (or $\sum_{k\in\mathbb{Z}} (-1)^k \hat{\gamma}(w+2k\pi)$) also converge pointwise everywhere when

$\sum_{k\in\mathbb{Z}} \hat{\gamma}(w+4k\pi)$ and $\sum_{k\in\mathbb{Z}} \hat{\gamma}(k/2)e^{-iwk/2}$ converge in (A2) (or in (A3)). Hence, we have

$$\begin{cases} \sum_{k\in\mathbb{Z}} \gamma(k)e^{-iwk} = \sum_{k\in\mathbb{Z}} \hat{\gamma}(w+2k\pi) \\ e^{-iw/2}\sum_{k\in\mathbb{Z}} \gamma(k+1/2)e^{-iwk} = \sum_{k\in\mathbb{Z}} (-1)^k \hat{\gamma}(w+2k\pi) \end{cases}. \tag{A4}$$

Since $\hat{S}^\gamma(w) = \hat{\gamma}(w)/\sum_{k\in\mathbb{Z}} \gamma(k)e^{-iwk}$, it follows from the first identity in (A4) that

$$\hat{S}^\gamma(w) = \frac{\hat{\gamma}(w)}{\sum_{k\in\mathbb{Z}} \hat{\gamma}(w+2k\pi)} \tag{A5}$$

and

$$\hat{S}^\gamma(w+4k\pi) = \frac{\hat{\gamma}(w+4k\pi)}{\sum_{n\in\mathbb{Z}} \hat{\gamma}(w+2n\pi)}, \; k\in\mathbb{Z}. \tag{A6}$$

Since the series $\sum_{k\in\mathbb{Z}} \hat{\gamma}(w+4k\pi)$ and $\sum_{k\in\mathbb{Z}} \hat{\gamma}(w+2k\pi)$ converge pointwise everywhere, equation (A6) indicates that $\sum_{k\in\mathbb{Z}} \hat{S}^\gamma(w+4k\pi) = \sum_{k\in\mathbb{Z}} \hat{\gamma}(w+4k\pi)/\sum_{n\in\mathbb{Z}} \hat{\gamma}(w+2n\pi)$ also converges pointwise everywhere.

Now, we show that $\frac{1}{2}\sum_{k\in\mathbb{Z}} S^\gamma(k/2)e^{-iwk/2}$ converges pointwise everywhere.

Since $\{S^\gamma(x-k)\}_{k\in\mathbb{Z}}$ and $\{\gamma(x-k)\}_{k\in\mathbb{Z}}$ are respectively an interpolation basis and a Riesz basis in the space $U_0(\gamma)$, there exists a unique set of coefficients $\{a_k\} \in l^2(\mathbb{R})$ such that

$$S^\gamma(x) = \sum_{k\in\mathbb{Z}} a_k \gamma(x-k). \tag{A7}$$

The Fourier transform of (A7) yields

$$\hat{S}^\gamma(w) = \sum_{k\in\mathbb{Z}} a_k e^{-iwk} \hat{\gamma}(w). \tag{A8}$$

The uniqueness of (A7) implies that the element $\sum_{k\in\mathbb{Z}} a_k e^{-iwk}$ of $L^2[-\pi,\pi]$ is unique in satisfying (A8). Observe that the function $1/\sum_{k\in\mathbb{Z}} \hat{\gamma}(w+2k\pi)$ is an element of $L^2[-\pi,\pi]$ if there exists an interpolation basis $\{S^\gamma(x-k)\}_k$ in $U_0(\gamma)$ (see Theorem 1 in [8]). Hence, it follows from (A5) and (A8) that

$$\sum_{k\in\mathbb{Z}} a_k e^{-iwk} = 1/\sum_{k\in\mathbb{Z}} \hat{\gamma}(w+2k\pi), \tag{A9}$$

Since $1/\sum_{k\in\mathbb{Z}} \hat{\gamma}(w+2k\pi)$ converges pointwise everywhere, equation (A9) indicates that $\sum_{k\in\mathbb{Z}} a_k e^{-iwk}$ does also.

It follows from (A7) that

$$S^\gamma(n/2) = \sum_{k\in\mathbb{Z}} a_k \gamma(n/2-k) \tag{A10}$$

By taking the Fourier transform on both sides in (A10), we obtain

$$\sum_{n\in\mathbb{Z}} S^\gamma(n/2)e^{-iwn/2} = 2\sum_{k\in\mathbb{Z}} a_k e^{-iwk} \sum_{n\in\mathbb{Z}} \hat{\gamma}(w+4n\pi) \tag{A11}$$

It follows from (A9) and (A11) that $\sum_n S^\gamma(n/2)e^{-iwn/2}$ converges pointwise everywhere since $\sum_{k\in\mathbb{Z}} a_k e^{-iwk}$ and $\sum_{n\in\mathbb{Z}} \hat{\gamma}(w+4n\pi)$ both do.

Since both $\sum_{n\in\mathbb{Z}} S^\gamma(n/2)e^{-iwn/2}$ and $\sum_{k\in\mathbb{Z}} \hat{S}^\gamma(w+4k\pi)$ converge pointwise everywhere, inserting (A9) into (A11) yields

$$\frac{1}{2}\sum_{k\in\mathbb{Z}} S^\gamma(k/2)z^k = \sum_{k\in\mathbb{Z}} \hat{\gamma}(w+4k\pi)/\sum_{k\in\mathbb{Z}} \hat{\gamma}(w+2n\pi) = \sum_{k\in\mathbb{Z}} \hat{S}^\gamma(w+4k\pi), \tag{A12}$$

which implies that the Poisson summation formula holds for $\{S^\gamma(k/2)\}_{k\in\mathbb{Z}}$.

**Part 2: Proof that** $P_s(z) = \frac{1}{2}\sum_{k\in\mathbb{Z}} S^\phi(k/2)z^k = \sum_{k\in\mathbb{Z}} \hat{S}^\phi(w+4k\pi)$ **with** $z = e^{-iw/2}$ **is the two-scale symbol corresponding to the scaling function** $S^\phi(x)$ **and**

$$P_s(z) = \frac{P_\phi(z)\sum_{k\in\mathbb{Z}}\hat{\phi}(w/2+2k\pi)}{\sum_{n\in\mathbb{Z}}\hat{\phi}(w+2n\pi)}.$$

Since $V_0 \subset V_1$ and $\{S^\phi(2x-k)\}_k$ is an interpolation basis with interpolation points located at $k/2$, $k \in \mathbb{Z}$, we have

$$S^\phi(x) = \sum_{k=-\infty}^{+\infty} S^\phi(k/2) S^\phi(2x-k). \tag{A13}$$

Since we have assumed that both $\sum_{k\in\mathbb{Z}} \phi(k/2)e^{-iwk/2}$ and $\sum_{k\in\mathbb{Z}} \hat{\phi}(w+4k\pi)$ in (3) converge pointwise everywhere, it follows from Part 1 of this proof that $\sum_{k\in\mathbb{Z}} S^\phi(k/2)z^k$ and $\sum_{k\in\mathbb{Z}} \hat{S}^\phi(w+4k\pi)$ do also. Hence, it follows from (A12) and (A13) that

$$\hat{S}^\phi(w) = P_s(z)\hat{S}^\phi(w/2) = \sum_{k=-\infty}^{+\infty} S^\phi(k/2)z^k \hat{S}(w/2)$$
$$= \sum_{k=-\infty}^{+\infty} \hat{S}^\phi(w+4k\pi)\hat{S}^\phi(w/2) \qquad , z = e^{-iw/2} \tag{A14}$$

Equation (A14) implies that

$$P_s(z) = \sum_{k=-\infty}^{+\infty} S^\phi(k/2)z^k = \sum_{k=-\infty}^{+\infty} \hat{S}^\phi(w+4k\pi). \tag{A15}$$

Applying (10) to (A15) gives

$$P_s(z) = \frac{\sum_{k\in\mathbb{Z}}\hat{\phi}(w+4k\pi)}{\sum_{n\in\mathbb{Z}}\hat{\phi}(w+2n\pi)} \tag{A16}$$

On the other hand, it follows from (8) that

$$\hat{\phi}(w+4k\pi) = P_\phi(e^{-i(w+4k\pi)/2})\hat{\phi}(w/2+2k\pi) = P_\phi(z)\hat{\phi}(w/2+2k\pi), \ k \in \mathbb{Z} \tag{A17}$$

Inserting (A17) into (A16) yields

$$P_s(z) = \frac{P_\phi(z)\sum_{k\in\mathbb{Z}}\hat{\phi}(w/2+2k\pi)}{\sum_{n\in\mathbb{Z}}\hat{\phi}(w+2n\pi)}. \tag{A18}$$

Equations (A14), (A15) and (A18) respectively imply (13), (12) and (14). □

### Appendix B : PROOF OF LEMMA 1

**Proof.** In (24), the change of variable $y = 2^j x$ gives

$$\sum_{k\in\mathbb{Z}} b_k^0 \tilde{S}^\phi(y-k) = \sum_{k\in\mathbb{Z}} c_k^0 \tilde{S}^\phi(2y-2k) + \sum_{k\in\mathbb{Z}} d_k^0 \tilde{S}^\phi(2y-2k-1). \tag{B1}$$

Obviously, equation (24) holds if and only if (B1) does.

**First, we show that**

$$0 < \left|\sum_k \tilde{p}_{2k}^s e^{-iwk}\right|^2 < +\infty, \text{ with } w \in [-\pi, \pi]$$

**for the coefficients $\{\tilde{p}_k^s\}_{k \in \mathbb{Z}}$ in (22).**

Since $P_{\tilde{S}^\phi}(z)$ is a two-scale symbol corresponding to the scaling function $\tilde{S}^\phi(x)$, it follows from (14) that

$$P_s(z) = \frac{P_{\tilde{S}^\phi}(z) \sum_{k \in \mathbb{Z}} \hat{\tilde{S}}^\phi(w/2 + 2k\pi)}{\sum_{n \in \mathbb{Z}} \hat{\tilde{S}}^\phi(w + 2n\pi)} \tag{B2}$$

Since $\{\tilde{S}^\phi(x-k)\}_{k \in \mathbb{Z}}$ forms the dual basis of $\{S^\phi(x-k)\}_{k \in \mathbb{Z}}$, we have

$$\hat{\tilde{S}}^\phi(w) = \frac{\hat{S}^\phi(w)}{\sum_{k=-\infty}^{+\infty} \left|\hat{S}^\phi(w + 2k\pi)\right|^2} \tag{B3}$$

(see equation (5.6.8) in [20]). Inserting (B3) into (B2) yields

$$P_s(z) = \frac{P_{\tilde{S}^\phi}(z) \sum_{k \in \mathbb{Z}} \hat{S}^\phi(w/2 + 2k\pi)}{\sum_{k \in \mathbb{Z}} \hat{S}^\phi(w + 2k\pi)} \times \frac{\sum_{k \in \mathbb{Z}} \left|\hat{S}^\phi(w + 2k\pi)\right|^2}{\sum_{k \in \mathbb{Z}} \left|\hat{S}^\phi(w/2 + 2k\pi)\right|^2}. \tag{B4}$$

It follows from (10) that

$$\sum_{k \in \mathbb{Z}} \hat{S}^\phi(w + 2k\pi) = 1. \tag{B5}$$

Inserting (B5) and (16) into (B4) yields

$$P_s(z) = P_{\tilde{S}^\phi}(z) \frac{\sum_{k \in \mathbb{Z}} \left|\hat{S}^\phi(w + 2k\pi)\right|^2}{\sum_{k \in \mathbb{Z}} \left|\hat{S}^\phi(w/2 + 2k\pi)\right|^2} = P_{\tilde{S}^\phi}(z) \frac{E_s(z^2)}{E_s(z)}. \tag{B6}$$

Equation (B6) indicates that

$$\sum_k \tilde{p}_{2k}^s e^{-iwk} = \frac{1}{2}\left(P_{\tilde{S}^\phi}(-z) + P_{\tilde{S}^\phi}(z)\right) = \frac{P_s(-z)E_s(-z) + P_s(z)E_s(z)}{2E_s(z^2)} \tag{B7}$$

It follows from (B7) that

$$\left|\sum_k \tilde{p}_{2k}^s e^{-iwk}\right|^2 = \left|\overline{\sum_k \tilde{p}_{2k}^s e^{-iwk}}\right|^2 = \left|\frac{P_s(-z)E_s(-z) + P_s(z)E_s(z)}{2E_s(z^2)}\right|^2 = \frac{\left|\overline{P_s(-z)}E_s(-z) + \overline{P_s(z)}E_s(z)\right|^2}{4\left|E_s(z^2)\right|^2} \tag{B8}$$

Since $\{S^\phi(x-k)\}_{k \in \mathbb{Z}}$ forms a Riesz basis of $V_0$, there exists two constants $A_{Es}$ and $B_{Es}$ with $0 < A_{Es} \leq B_{Es} < +\infty$ such that

$$0 < A_{Es} \leq E_s(z^2) \leq B_{Es} < +\infty \text{ (Theorem 9 in [19])}. \tag{B9}$$

It follows from (23), (B8) and (B9) that

$$0 < \frac{A_s}{4B_{Es}^2} \leq \left|\sum_k \tilde{p}_{2k}^s e^{-iwk}\right|^2 \leq \frac{B_s}{4A_{Es}^2} < +\infty \qquad (B10)$$

Equation (B10) indicate that $0 < \left|\sum_k \tilde{p}_{2k}^s e^{-iwk}\right|^2 < +\infty$ for every $w \in [-\pi, \pi]$.

**Second, we show that there always exists a unique set of coefficients $\{b_k^0\}_{k\in\mathbb{Z}}$ and $\{d_k^0\}_{k\in\mathbb{Z}}$ for which (B1) holds.**

Assume that $\{c_k^0\}_{k\in\mathbb{Z}}$ is any element of $l^2(\mathbb{R})$. Since $0 < \left|\sum_k \tilde{p}_{2k}^s e^{-iwk}\right|^2 < +\infty$ by (B10), the series $2\sum_{k\in\mathbb{Z}} c_k^0 e^{-iwk} / \sum_{k\in\mathbb{Z}} \tilde{p}_{2k}^s e^{-iwk}$ is an element of $L^2[-\pi,\pi]$ for any given $\{c_k^0\}_{k\in\mathbb{Z}} \in l^2(\mathbb{R})$.

Since $2\sum_{k\in\mathbb{Z}} c_k^0 e^{-iwk} / \sum_{k\in\mathbb{Z}} \tilde{p}_{2k}^s e^{-iwk} \in L^2[-\pi,\pi]$, there always exists a unique set of coefficients $\{b_k^0\}_k \in L^2(\mathbb{R})$ such that

$$\sum_{k\in\mathbb{Z}} b_k^0 e^{-iwk} = 2\sum_{k\in\mathbb{Z}} c_k^0 e^{-iwk} / \sum_{k\in\mathbb{Z}} \tilde{p}_{2k}^s e^{-iwk}. \qquad (B11)$$

On the other hand, since $\{b_k^0\}_k \in l^2(\mathbb{R})$ in (B11), we have that $\frac{1}{2}\sum_{k\in\mathbb{Z}} b_k^0 e^{-iwk} \sum_{k\in\mathbb{Z}} \tilde{p}_{2k+1}^s e^{-iwk}$ is also an element of $L^2[-\pi,\pi]$. Hence, there always exists a unique set of coefficients $\{d_k^0\}_k \in l^2(\mathbb{R})$ such that

$$\frac{1}{2}\sum_{k\in\mathbb{Z}} b_k^0 e^{-iwk} \sum_{k\in\mathbb{Z}} \tilde{p}_{2k+1}^s e^{-iwk} = \sum_{k\in\mathbb{Z}} d_k^0 e^{-iwk}. \qquad (B12)$$

Equations (B11) and (B12) indicate that there always exist two unique sets of coefficients $\{b_k^0\}_{k\in\mathbb{Z}}$ and $\{d_k^0\}_{k\in\mathbb{Z}}$ in $l^2(\mathbb{R})$ such that

$$\begin{aligned}\sum_{k\in\mathbb{Z}} c_k^0 e^{-iwk} + e^{-iw/2}\sum_{k\in\mathbb{Z}} d_k^0 e^{-iwk} &= \frac{1}{2}\sum_{k\in\mathbb{Z}} b_k^0 e^{-iwk} \sum_{k\in\mathbb{Z}} \tilde{p}_{2k}^s e^{-iwk} \\ &+ \frac{1}{2}e^{-iw/2}\sum_{k\in\mathbb{Z}} b_k^0 e^{-iwk} \sum_{k\in\mathbb{Z}} \tilde{p}_{2k+1}^s e^{-iwk} \\ &= \sum_{k\in\mathbb{Z}} b_k^0 e^{-iwk} \left(\sum_{k\in\mathbb{Z}} \tilde{p}_{2k}^s e^{-iwk} + e^{-iw/2}\sum_{k\in\mathbb{Z}} \tilde{p}_{2k+1}^s e^{-iwk}\right)/2 \\ &= \frac{1}{2}\sum_{k\in\mathbb{Z}} b_k^0 e^{-iwk} \sum_{n\in\mathbb{Z}} \tilde{p}_n^s e^{-iwn/2}\end{aligned} \qquad (B13)$$

for any given set of $\{c_k^0\}_{k\in\mathbb{Z}} \in l^2(\mathbb{R})$.

Multiplying (B13) by $\hat{\tilde{S}}^\phi(w/2)$ and inserting (22) yield

$$\begin{aligned}\left(\sum_{k\in\mathbb{Z}} c_k^0 e^{-iwk} + e^{-iw/2}\sum_{k\in\mathbb{Z}} d_k^0 e^{-iwk}\right)\hat{\tilde{S}}^\phi(w/2) &= \frac{1}{2}\sum_{k\in\mathbb{Z}} b_k^0 e^{-iwk} \sum_{n\in\mathbb{Z}} \tilde{p}_n^s e^{-iwn/2} \hat{\tilde{S}}^\phi(w/2) \\ &= \sum_{k\in\mathbb{Z}} b_k^0 e^{-iwk} \hat{\tilde{S}}^\phi(w)\end{aligned}. \qquad (B14)$$

for any $\{c_k^0\}_k \in l^2(\mathbb{R})$.

Taking the inverse Fourier transform on both sides in (B14), we obtain

$$\sum_k c_k^0 \tilde{S}^\phi(2x-2k) + \sum_k d_k^0 \tilde{S}^\phi(2x-2k-1) = \sum_k b_k^0 \tilde{S}^\phi(x-k). \tag{B15}$$

Since the coefficients $\{b_k^0\}_k$ and $\{d_k^0\}_k$ for which (B13) holds are unique, it follows from (B15) that, for any given $\{c_k^0\}_k \in l^2(\mathbb{R})$, there always exist two unique sets of $\{b_k^0\}_{k\in\mathbb{Z}}$, $\{d_k^0\}_{k\in\mathbb{Z}} \in l^2(\mathbb{R})$ for which (B1) holds. Hence, for a given set of coefficients $\{c_k^0\}_k \in l^2(\mathbb{R})$, there always exists two unique sets of coefficients $\{b_k^0\}_{k\in\mathbb{Z}}$, $\{d_k^0\}_{k\in\mathbb{Z}} \in l^2(\mathbb{R})$ for which equation (24) holds.□

### Appendix C : PROOF OF THEOREM 1

**Proof.** Let $f_s(x) \in V_j$. From Section 2, we know that if there exists an interpolation basis $\{S^\phi(x-k)\}_k$ in $V_0$, then there must exist an interpolation basis $\{S^\phi(2^j x - k)\}_k$ in every $V_j$. Hence, the function $f_s(x) \in V_j$ can be written as

$$f_s(x) = \sum_{k\in\mathbb{Z}} f_s(\frac{k}{2^j}) S^\phi_{j,k}(x) = \sum_{k\in\mathbb{Z}} \left( f_s(\frac{2k+1}{2^j}) S^\phi_{j,2k+1}(x) + f_s(\frac{2k}{2^j}) S^\phi_{j,2k}(x) \right) \tag{C1}$$

where $S^\phi_{j,k}(x) = S^\phi(2^j x - k)$.

Substituting (25) into (C1) yields

$$f_s(x) = \sum_{k\in\mathbb{Z}} f_s(\frac{k}{2^j}) S^\phi_{j,k}(x) = \sum_{k\in\mathbb{Z}} \left( F_k S^\phi_{j,2k+1}(x) + f_s(\frac{2k}{2^j}) S^\phi_{j,2k}(x) \right). \tag{C2}$$

Inserting (C2) into $\sum_{k\in\mathbb{Z}} \langle f_s, S^\phi_{j-1,k} \rangle \tilde{S}^\phi_{j-1,k}$ yields

$$\sum_{n\in\mathbb{Z}} \langle f_s, S^\phi_{j-1,n} \rangle \tilde{S}^\phi_{j-1,n} = \sum_{n\in\mathbb{Z}}\sum_{k\in\mathbb{Z}} F_k \langle S^\phi_{j,2k+1}(x), S^\phi_{j-1,n} \rangle \tilde{S}^\phi_{j-1,n} + \sum_{n\in\mathbb{Z}}\sum_{k\in\mathbb{Z}} f_s(\frac{2k}{2^j}) \langle S^\phi_{j,2k}(x), S^\phi_{j-1,n} \rangle \tilde{S}^\phi_{j-1,n} \tag{C3}$$

where $\{\tilde{S}^\phi_{j-1,k}(x) = \tilde{S}^\phi(2^{j-1} x - k)\}_k$ is the dual basis of $\{S^\phi_{j-1,k}(x) = S^\phi(2^{j-1} x - k)\}_k$.

Obviously, $\sum_{n\in\mathbb{Z}} \langle f_s, S^\phi_{j-1,n} \rangle \tilde{S}^\phi_{j-1,n}$ in (C3) is the projection of $f_s(x)$ on $V_{j-1}$. Hence, $f_s(x)$ is an element of $W_{j-1}$ if and only if $\sum_{n\in\mathbb{Z}} \langle f_s, S^\phi_{j-1,n} \rangle \tilde{S}^\phi_{j-1,n} = 0$ in (C3).

Now, we show there exists a unique set of $\{f_s(2k/2^j)\}_{k\in\mathbb{Z}}$ in (C3) such that

$$\sum_{n\in\mathbb{Z}} \langle f_s, S^\phi_{j-1,n} \rangle \tilde{S}^\phi_{j-1,n} = \sum_{n\in\mathbb{Z}}\sum_{k\in\mathbb{Z}} F_k \langle S^\phi_{j,2k+1}(x), S^\phi_{j-1,n} \rangle \tilde{S}^\phi_{j-1,n} + \sum_{n\in\mathbb{Z}}\sum_{k\in\mathbb{Z}} f_s(\frac{2k}{2^j}) \langle S^\phi_{j,2k}(x), S^\phi_{j-1,n} \rangle \tilde{S}^\phi_{j-1,n} = 0 \tag{C4}$$

for a given set of $\{F_k\}_k$, i.e. the sequence $\{F_k\}_k$ in (C2) determines a unique function $f_s(x)$ whose projection $\sum_{n\in\mathbb{Z}} \langle f_s, S^\phi_{j-1,n}\rangle \tilde{S}^\phi_{j-1,n}$ on $V_{j-1}$ vanishes.

Letting

$$a_l = \langle S^\phi(2^j x), S^\phi(2^j x - l)\rangle = \frac{1}{2^j}\langle S^\phi(x), S^\phi(x-l)\rangle, \tag{C5}$$

then

$$S^\phi(2^j x - k) = \sum_l \langle S^\phi(2^j x), S^\phi(2^j x - l)\rangle \tilde{S}^\phi(2^j x - l - k) = \sum_l a_l \tilde{S}^\phi(2^j x - l - k). \tag{C6}$$

Substituting (C6) into (C2) yields that

$$f_s(x) = \sum_{k\in\mathbb{Z}} F_k \sum_{n\in\mathbb{Z}} a_n \tilde{S}^\phi(2^j x - n - 2k - 1) + \sum_{k\in\mathbb{Z}} f_s(\frac{2k}{2^j}) \sum_{n\in\mathbb{Z}} a_n \tilde{S}^\phi(2^j x - n - 2k) \tag{C7}$$

It follows from (C4) and (C7) that

$$\langle f_s, S^\phi_{j-1,l}\rangle = \sum_{k\in\mathbb{Z}} F_k \sum_{n\in\mathbb{Z}} a_n \langle \tilde{S}^\phi_{j,n+2k+1}, S^\phi_{j-1,l}\rangle + \sum_{k\in\mathbb{Z}} f_s(\frac{2k}{2^j}) \sum_{n\in\mathbb{Z}} a_n \langle \tilde{S}^\phi_{j,n+2k}, S^\phi_{j-1,l}\rangle = 0. \tag{C8}$$

Due to the fact that $\{\tilde{S}^\phi_{j,k}(x)\}_k$ is the dual basis of the interpolation basis $\{S^\phi_{j,k}(x)\}_k$ in $V_j$, we have

$$\begin{cases} \langle \tilde{S}^\phi_{j,n+2k+1}, S^\phi_{j-1,l}\rangle = \int_{-\infty}^{+\infty} \tilde{S}^\phi_{j,n+2k+1}(x)\overline{S}^\phi_{j-1,l}(x)dx = \overline{S}^\phi(2^{j-1}\frac{n+2k+1}{2^j} - l) = \overline{S}^\phi(\frac{n+2k-2l+1}{2}) \\ \langle \tilde{S}^\phi_{j,n+2k}, S^\phi_{j-1,l}\rangle = \int_{-\infty}^{+\infty} \tilde{S}^\phi_{j,n+2k}(x)\overline{S}^\phi_{j-1,l}(x)dx = \overline{S}^\phi(2^{j-1}\frac{n+2k}{2^j} - l) = \overline{S}^\phi(\frac{n+2k-2l}{2}) \end{cases} \tag{C9}$$

where $\overline{S}^\phi(x)$ is the complex conjugate of $S^\phi(x)$.

Hence, equation (C8) can be rewritten as

$$\langle f_s, S^\phi_{j-1,l}\rangle = \sum_{k\in\mathbb{Z}} F_k \sum_{n\in\mathbb{Z}} a_n \overline{S}^\phi(\frac{n+2k-2l+1}{2}) + \sum_{k\in\mathbb{Z}} f_s(\frac{2k}{2^j}) \sum_{n\in\mathbb{Z}} a_n \overline{S}^\phi(\frac{n+2k-2l}{2}) = 0. \tag{C10}$$

The Fourier transform of (C10) yields

$$\sum_{l\in\mathbb{Z}} \langle f_s, S^\phi_{j-1,l}\rangle e^{-iwl} = A_f + B_f = \sum_{k\in\mathbb{Z}} F_k \sum_{n\in\mathbb{Z}} a_n \sum_{l\in\mathbb{Z}} \overline{S}^\phi(\frac{n+2k-2l+1}{2})e^{-iwl}$$
$$+ \sum_{k\in\mathbb{Z}} f_s(\frac{2k}{2^j}) \sum_{n\in\mathbb{Z}} a_n \sum_{l\in\mathbb{Z}} \overline{S}^\phi(\frac{n+2k-2l}{2})e^{-iwl} = 0 \tag{C11}$$

where $A_f = \sum_{k\in\mathbb{Z}} F_k \sum_{n\in\mathbb{Z}} a_n \sum_{l\in\mathbb{Z}} \overline{S}^\phi(\frac{n+2k-2l+1}{2})e^{-iwl}$ and $B_f = \sum_{k\in\mathbb{Z}} f_s(\frac{2k}{2^j}) \sum_{n\in\mathbb{Z}} a_n \sum_{l\in\mathbb{Z}} \overline{S}^\phi(\frac{n+2k-2l}{2})e^{-iwl}$.

Similarly to $\phi(x)$ in (4), since the Poisson summation formula holds for sequence $\{S^\phi(k/2)\}_{k\in\mathbb{Z}}$ by Proposition 1, the Poisson summation formula also holds for sequence $\{S^\phi(k+1/2)\}_{k\in\mathbb{Z}}$. Hence, we have

$$\sum_{l\in\mathbb{Z}}\overline{S^{\phi}(\frac{n+2k-2l+1}{2})}e^{-iwl} = \overline{\sum_{l\in\mathbb{Z}}S^{\phi}(\frac{n+1}{2}+k-l)e^{iwl}} = \overline{\sum_{l\in\mathbb{Z}}S^{\phi}(\frac{n+1}{2}+k+l)e^{-iwl}}$$
$$= \overline{\sum_{l\in\mathbb{Z}}\hat{S}^{\phi}(w+2l\pi)e^{i(w+2l\pi)(n+1+2k)/2}} \quad . \tag{C12}$$

where $\hat{S}^{\phi}(w)$ is the Fourier transform of $S^{\phi}(x)$.

Inserting (C12) in $A_f$ yields

$$A_f = \sum_{k\in\mathbb{Z}}F_k e^{-ikw}\sum_{n\in\mathbb{Z}}a_n e^{-\frac{n+1}{2}wi}\overline{\sum_{l\in\mathbb{Z}}\hat{S}^{\phi}(w+2l\pi)e^{i(n+1)l\pi}} \quad . \tag{C13}$$

Cutting the summation (C13) into two parts, one for $n = 2m$ and the other for $n = 2m-1$, yields

$$A_f = \sum_{k\in\mathbb{Z}}F_k e^{-ikw}\left(\sum_{m\in\mathbb{Z}}a_{2m}e^{-mwi}\overline{e^{iw/2}\sum_{l\in\mathbb{Z}}\hat{S}^{\phi}(w+2l\pi)(-1)^l} + \sum_{m\in\mathbb{Z}}a_{2m-1}e^{-mwi}\overline{\sum_{l\in\mathbb{Z}}\hat{S}^{\phi}(w+2l\pi)}\right). \tag{C14}$$

Since $\{S^{\phi}(x-k)\}_k$ is the interpolation basis of $V_0$, it follows from (10) that

$$\sum_{l\in\mathbb{Z}}\hat{S}^{\phi}(w+2l\pi) = 1 \tag{C15}$$

Substituting (C15) into (C14) yields

$$A_f = \sum_{k\in\mathbb{Z}}F_k e^{-ikw}\left(\sum_{m\in\mathbb{Z}}a_{2m-1}e^{-mwi} + \sum_{m\in\mathbb{Z}}a_{2m}e^{-mwi}\overline{e^{-iw/2}\sum_{l\in\mathbb{Z}}\hat{S}^{\phi}(w+2l\pi)(-1)^l}\right). \tag{C16}$$

It follows from (C5) that

$$\sum_{m\in\mathbb{Z}}a_{2m}e^{-iwm} = \frac{1}{2^j}\langle S^{\phi}(x), S^{\phi}(x-2m)\rangle = \frac{1}{2^j}\sum_{m\in\mathbb{Z}}\int_{-\infty}^{+\infty}\overline{S^{\phi}}(x)S^{\phi}(x-2m)dx \times e^{-iwm}$$
$$= \frac{1}{2^{j+1}}\frac{1}{2\pi}\sum_{m\in\mathbb{Z}}\int_{-\infty}^{+\infty}|\hat{S}^{\phi}(u/2)|^2 e^{-ium}du \times e^{-iwm} \quad . \tag{C17}$$

In sense of distributions (Theorem 2.4 in [6]), equation (C17) can be rewritten as

$$\sum_{m\in\mathbb{Z}}a_{2m}e^{-iwm} = \frac{1}{2^{j+1}}\frac{1}{2\pi}\int_{-\infty}^{+\infty}|\hat{S}^{\phi}(u/2)|^2 \sum_{m\in\mathbb{Z}}e^{-ium}\times e^{-iwm}du$$
$$= \frac{1}{2^{j+1}}\int_{-\infty}^{+\infty}|\hat{S}^{\phi}(u/2)|^2 \sum_{m\in\mathbb{Z}}\delta(u+w-2m\pi)du \quad . \tag{C18}$$

Observe that $\{S^{\phi}(x-k)\}_{k\in\mathbb{Z}}$ is a Riesz basis of $V_0$, so

$$\sum_{k\in\mathbb{Z}}|\hat{S}^{\phi}(\frac{w+2k\pi}{2})|^2 = \sum_{k\in\mathbb{Z}}|\hat{S}^{\phi}(w/2+2k\pi)|^2 + \sum_{k\in\mathbb{Z}}|\hat{S}^{\phi}(w/2+\pi+2k\pi)|^2$$

converges by the frame theory [19]. Since $|\hat{S}^{\phi}(u/2)|^2$ is an even function, it follows from (C18) that

$$\sum_{m\in\mathbb{Z}}a_{2m}e^{-mwi} = \frac{1}{2^{j+1}}\sum_{k\in\mathbb{Z}}\left|\hat{S}^{\phi}(\frac{w+2k\pi}{2})\right|^2 . \tag{C19}$$

Similarly, we have

$$\sum_{m\in\mathbb{Z}} a_{2m-1} e^{-mwi} = \frac{1}{2^{j+1}} e^{-\frac{iw}{2}} \sum_{k\in\mathbb{Z}} (-1)^k \left| \hat{S}^\phi\left(\frac{w+2k\pi}{2}\right) \right|^2. \tag{C20}$$

Inserting (C19) and (C20) into (C16) yields

$$A_f = \frac{e^{-iw/2}}{2^{j+1}} \sum_{k\in\mathbb{Z}} F_k e^{-ikw} \left( \sum_{m\in\mathbb{Z}} (-1)^m \left| \hat{S}^\phi\left(\frac{w+2m\pi}{2}\right) \right|^2 \right.$$
$$\left. + \sum_{m\in\mathbb{Z}} \left| \hat{S}^\phi\left(\frac{w+2m\pi}{2}\right) \right|^2 \overline{\left( \sum_{l\in\mathbb{Z}} \hat{S}^\phi(w+2l\pi)(-1)^l \right)} \right). \tag{C21}$$

Cutting the summation (C21) into two parts, one for $m=2n$ and the other for $m=2n-1$, yields

$$A_f = \frac{e^{-iw/2}}{2^{j+1}} \sum_{k\in\mathbb{Z}} F_k e^{-ikw} \times \left( \sum_{n\in\mathbb{Z}} \left| \hat{S}^\phi(w/2 + 2\pi n) \right|^2 \overline{\left( 1 + \sum_{l\in\mathbb{Z}} (-1)^l \hat{S}^\phi(w+2l\pi) \right)} \right.$$
$$\left. - \sum_{n\in\mathbb{Z}} \left| \hat{S}^\phi(w/2 + 2\pi n + \pi) \right|^2 \overline{\left( 1 - \sum_{l\in\mathbb{Z}} (-1)^l \hat{S}^\phi(w+2l\pi) \right)} \right) . \tag{C22}$$

Since it follows from (10) that $\sum_{l\in\mathbb{Z}} \hat{S}^\phi(w+2l\pi) = \sum_{l\in\mathbb{Z}} \hat{S}^\phi(w+4l\pi) + \sum_{l\in\mathbb{Z}} \hat{S}^\phi(w+4l\pi+2\pi) = 1$, equation (C22) indicates that

$$A_f = \frac{e^{-iw/2}}{2^j} \sum_{k\in\mathbb{Z}} F_k e^{-ikw} \times \left( \sum_{n\in\mathbb{Z}} \left| \hat{S}^\phi(w/2 + 2\pi n) \right|^2 \overline{\sum_{l\in\mathbb{Z}} \hat{S}^\phi(w+4l\pi)} \right.$$
$$\left. - \sum_{n\in\mathbb{Z}} \left| \hat{S}^\phi(w/2 + 2\pi n + \pi) \right|^2 \overline{\sum_{l\in\mathbb{Z}} \hat{S}^\phi(w+4l\pi+2\pi)} \right) . \tag{C23}$$

Similarly,

$$B_f = \frac{1}{2^j} \sum_{k\in\mathbb{Z}} f_s\left(\frac{2k}{2^j}\right) e^{-ikw} \times \left( \sum_{n\in\mathbb{Z}} \left| \hat{S}^\phi(w/2 + 2\pi n) \right|^2 \overline{\left( \sum_{l\in\mathbb{Z}} \hat{S}^\phi(w+4l\pi) \right)} \right.$$
$$\left. + \sum_{n\in\mathbb{Z}} \left| \hat{S}^\phi(w/2 + 2\pi n + \pi) \right|^2 \overline{\left( \sum_{l\in\mathbb{Z}} \hat{S}^\phi(w+4l\pi+2\pi) \right)} \right) . \tag{C24}$$

Let

$$\tilde{Q}_s(z) = -z E_s(-z) \overline{P_s(-z)} \text{ with } z = e^{-\frac{wi}{2}} \tag{C25}$$

and

$$\hat{\psi}_Q(w) = \tilde{Q}_s(z) \hat{S}^\phi(w/2), \tag{C26}$$

where $P_s(z)$ and $E_s(z)$ are respectively defined in (12) and (16).

Since $P_s(z)$ is a two-scale symbol corresponding to the scaling function $S^\phi(2^j x)$ by Proposition 1, the sequence $\{\psi_Q(x-k)\}_k$ forms a Riesz basis of $W_0$ (see equation (5.6.13) in [20]), where the Fourier transform of $\psi_Q(x)$ is $\hat{\psi}_Q(w)$.

It follows from (12) that

$$P_s(z) = \sum_{n \in \mathbb{Z}} \hat{S}^\phi(w - 4n\pi). \tag{C27}$$

Inserting (C27) and (16) into (C25) yields

$$\tilde{Q}_s(z) = -z \sum_{k \in \mathbb{Z}} \left| \hat{S}^\phi(\frac{w}{2} + \pi + 2\pi k) \right|^2 \overline{\sum_{n \in \mathbb{Z}} \hat{S}^\phi(w + 2\pi + 4n\pi)} \tag{C28}$$

and

$$\tilde{Q}_s(-z) = z \sum_{k \in \mathbb{Z}} \left| \hat{S}^\phi(\frac{w}{2} + 2\pi k) \right|^2 \overline{\sum_{n \in \mathbb{Z}} \hat{S}^\phi(w + 4n\pi)}. \tag{C29}$$

Inserting (C28) and (C29) into (C23) yields

$$A_f = \frac{1}{2^j} \sum_{k \in \mathbb{Z}} F_k e^{-ikw} \times \left( \tilde{Q}_s(-z) + \tilde{Q}_s(z) \right). \tag{C30}$$

Inserting (C28) and (C29) into (C24) yields

$$B_f = \frac{1}{2^j z} \sum_{k \in \mathbb{Z}} f_s(\frac{2k}{2^j}) e^{-ikw} \times \left( \tilde{Q}_s(-z) - \tilde{Q}_s(z) \right). \tag{C31}$$

Since $\sum_{k \in \mathbb{Z}} \left| \hat{S}^\phi(w/2 + 2\pi k) \right|^2 \in L^2[-\pi, \pi]$ by the frame theory [19] and since it follows from (11) and (12) that $\sum_{n \in \mathbb{Z}} \hat{S}^\phi(w + 4n\pi) \in L^2[-2\pi, 2\pi]$, $\tilde{Q}_s(z)$ can be written as

$$\tilde{Q}_s(z) = -z \sum_{k \in \mathbb{Z}} \left| \hat{S}^\phi(\frac{w}{2} + \pi + 2\pi k) \right|^2 \overline{\sum_{n \in \mathbb{Z}} \hat{S}^\phi(w + 2\pi + 4n\pi)} = \frac{1}{2} \sum_{k \in \mathbb{Z}} \tilde{q}_k z^k \tag{C32}$$

with $\{\tilde{q}_k\}_{k \in \mathbb{Z}} \in l^2(\mathbb{Z})$.

Then, applying (C32) to (C30) gives

$$A_f = \frac{1}{2^j} \sum_{k \in \mathbb{Z}} F_k e^{-ikw} \times \left( \tilde{Q}_s(-z) + \tilde{Q}_s(z) \right) = \frac{1}{2^j} \sum_{k \in \mathbb{Z}} F_k z^{2k} \times \sum_{l \in \mathbb{Z}} \tilde{q}_{2l} z^{2l}. \tag{C33}$$

Applying (C32) to (C31) gives

$$B_f = \frac{1}{2^j z} \sum_{k \in \mathbb{Z}} f_s(\frac{2k}{2^j}) e^{-ikw} \times \left( \tilde{Q}_s(-z) - \tilde{Q}_s(z) \right) = -\frac{1}{2^j} \sum_{k \in \mathbb{Z}} f_s(\frac{2k}{2^j}) z^{2k} \times \sum_{l \in \mathbb{Z}} \tilde{q}_{2l+1} z^{2l}. \tag{C34}$$

Substituting (C33) and (C34) into (C11) gives

$$\sum_{k\in\mathbb{Z}} f_s(\frac{2k}{2^j})z^{2k-1} \times \sum_{l\in\mathbb{Z}} \tilde{q}_{2l-1}z^{2l-1} - \sum_{k\in\mathbb{Z}} F_k z^{2k} \times \sum_{l\in\mathbb{Z}} \tilde{q}_{2l}z^{2l} = 0. \tag{C35}$$

Let $\tilde{\psi}_Q(x)$ be the dual of $\psi_Q(x)$ in $W_0$. Consider the "decomposition relation" of scaling function and wavelet, it follows from (12), (C32) and Theorem 5.16 in [20] that

$$\begin{cases} \hat{\tilde{S}}^\phi(\frac{w}{2}) = \sum_{k\in\mathbb{Z}} \left( S^\phi(\frac{2k}{2})z^{2k}\hat{\tilde{S}}^\phi(w) + \tilde{q}_{2k}z^{2k}\hat{\tilde{\psi}}_Q(w) \right) \\ \hat{\tilde{S}}^\phi(\frac{w}{2}) = \sum_{k\in\mathbb{Z}} \left( S^\phi(\frac{2k-1}{2})z^{2k-1}\hat{\tilde{S}}^\phi(w) + \tilde{q}_{2k-1}z^{2k-1}\hat{\tilde{\psi}}_Q(w) \right). \end{cases} \tag{C36}$$

where $\hat{\tilde{S}}^\phi(w)$ is the Fourier transform of $\tilde{S}^\phi(x)$ and $\hat{\tilde{\psi}}_Q(w)$ is the Fourier transform of $\tilde{\psi}_Q(x)$.

By multiplying the two identities in (C36) by $\sum_{k\in\mathbb{Z}}(-F_k)z^{2k}$ and $\sum_{k\in\mathbb{Z}} f_s(\frac{2k}{2^j})z^{2k-1}$ consecutively, we have

$$\begin{cases} -\sum_{k\in\mathbb{Z}} F_k z^{2k}\hat{\tilde{S}}^\phi(\frac{w}{2}) = -\sum_{k\in\mathbb{Z}} F_k z^{2k} \sum_{l\in\mathbb{Z}} S^\phi(\frac{2l}{2})z^{2l}\hat{\tilde{S}}^\phi(w) - \sum_{k\in\mathbb{Z}} F_k z^{2k} \sum_{l\in\mathbb{Z}} \tilde{q}_{2l}z^{2l}\hat{\tilde{\psi}}_Q(w) \\ \sum_{k\in\mathbb{Z}} f_s(\frac{2k}{2^j})z^{2k-1}\hat{\tilde{S}}^\phi(\frac{w}{2}) = \sum_{k\in\mathbb{Z}} f_s(\frac{2k}{2^j})z^{2k-1} \sum_{l\in\mathbb{Z}} S^\phi(\frac{2l-1}{2})z^{2l-1}\hat{\tilde{S}}^\phi(w) \\ \qquad + \sum_{k\in\mathbb{Z}} f_s(\frac{2k}{2^j})z^{2k-1} \sum_{l\in\mathbb{Z}} \tilde{q}_{2l-1}z^{2l-1}\hat{\tilde{\psi}}_Q(w) \end{cases} \tag{C37}$$

Then, it follows from (C37) that

$$-\sum_{k\in\mathbb{Z}} F_k z^{2k}\hat{\tilde{S}}^\phi(\frac{w}{2}) + \sum_{k\in\mathbb{Z}} f_s(\frac{2k}{2^j})z^{2k-1}\hat{\tilde{S}}^\phi(\frac{w}{2}) =$$
$$\left( -\sum_{k\in\mathbb{Z}} F_k z^{2k} \sum_{l\in\mathbb{Z}} S^\phi(\frac{2l}{2})z^{2l} + \sum_{k\in\mathbb{Z}} f_s(\frac{2k}{2^j})z^{2k-1} \sum_{l\in\mathbb{Z}} S^\phi(\frac{2l-1}{2})z^{2l-1} \right)\hat{\tilde{S}}^\phi(w)$$
$$+ \left( \sum_{k\in\mathbb{Z}} f_s(\frac{2k}{2^j})z^{2k-1} \sum_{l\in\mathbb{Z}} \tilde{q}_{2l-1}z^{2l-1} - \sum_{k\in\mathbb{Z}} F_k z^{2k} \sum_{l\in\mathbb{Z}} \tilde{q}_{2l}z^{2l} \right)\hat{\tilde{\psi}}_Q(w) \tag{C38}$$

Applying (C35) to (C38), we obtain

$$-\sum_{k\in\mathbb{Z}} F_k z^{2k}\hat{\tilde{S}}^\phi(\frac{w}{2}) + \sum_{k\in\mathbb{Z}} f_s(\frac{2k}{2^j})z^{2k-1}\hat{\tilde{S}}^\phi(\frac{w}{2}) =$$
$$\left( -\sum_{k\in\mathbb{Z}} F_k z^{2k} \sum_{l\in\mathbb{Z}} S^\phi(\frac{2l}{2})z^{2l} + \sum_{k\in\mathbb{Z}} f_s(\frac{2k}{2^j})z^{2k-1} \sum_{l\in\mathbb{Z}} S^\phi(\frac{2l-1}{2})z^{2l-1} \right)\hat{\tilde{S}}^\phi(w). \tag{C39}$$

By taking the inverse Fourier transform on both sides in (C39), we obtain

$$f_c(x) = -\sum_{k\in\mathbb{Z}} F_k \tilde{S}^\phi(2x-2k) + \sum_{k\in\mathbb{Z}} f_s(2k/2^j)\tilde{S}^\phi(2x-2k+1)$$
$$= \sum_{m\in\mathbb{Z}} \left( \sum_{n\in\mathbb{Z}} \left( f_s(\frac{2n}{2^j})S^\phi(\frac{2m-2n+1}{2}) - F_n S^\phi(\frac{2m-2n}{2}) \right) \times \tilde{S}^\phi(x-m) \right). \tag{C40}$$

Obviously, it follows from (C40) that $f_c(x)$ is an element of $V_0$. Comparing (C40) to (24), it is shown that $\{-F_k\}_{k\in\mathbb{Z}}$ and $\{f_s(2k/2^j)\}_{k\in\mathbb{Z}}$ in (C40) respectively corresponds to $\{c_k^0\}_{k\in\mathbb{Z}}$ and $\{d_k^0\}_{k\in\mathbb{Z}}$ in (24). Since equation (23) holds, it follows from Lemma 1 that the coefficients $\{f_s(2k/2^j)\}_{k \in \mathbb{Z}}$ for which (C40) holds is unique for any given sequence $\{-F_k\}_{k \in \mathbb{Z}}$. Hence, equation (C40) indicates that $\{F_k\}_{k \in \mathbb{Z}}$ determine uniquely $\{f_s(2k/2^j)\}_{k \in \mathbb{Z}}$ in (C4).

Since $\{f_s(2k/2^j)\}_{k \in \mathbb{Z}}$ and $\{F_k\}_{k \in \mathbb{Z}}$ determine a unique element of $V_j$ in (C2), equations (C4) and (C40) imply that the sequence $\{F_k\}_{k \in \mathbb{Z}}$ determines a unique element of $V_j$ whose projection on $V_{j-1}$ vanishes, i.e. $\{F_k\}_{k \in \mathbb{Z}}$ determines a unique element of $W_{j-1}$. This proves the statement in Theorem 1. □

**Appendix D : Proof of Proposition 3**

**Part 1: Proof that the inverse of $\tilde{T}_{W_{j-1}}$ is given by $f_s(x) = \sum_{k \in \mathbb{Z}} f_s(k/2^{j-1} + 1/2^j) S^\psi(2^{j-1}x - k)$.**

Obviously, it follows from (26) that $\tilde{W}_{j-1} \subset V_j$ (see equation (7.1) in [6]). Simultaneously, a translation

$$Tr_j(f(x)) = f(x - 1/2^j) \tag{D1}$$

is a bijective isometry of $L^2(\mathbb{R})$ onto $L^2(\mathbb{R})$, so the Hilbert spaces $\tilde{W}_{j-1}$ and $W_{j-1}$ are isometric in $V_j$.

Assuming that the sampling operator $T_{\tilde{W}_{j-1}}(f_s) = \{f_s(k/2^{j-1})\}_{k \in \mathbb{Z}}$ is a transformation of $\tilde{W}_{j-1}$ into $l^2(\mathbb{R})$, then it follows from (21) and (D1) that

$$\tilde{T}_{W_{j-1}} Tr_j(f_s) = \tilde{T}_{W_{j-1}}(f_s(x - 1/2^j)) = \{f_s(k/2^{j-1})\}_k = T_{\tilde{W}_{j-1}}(f_s), \tag{D2}$$

for any $f_s \in \tilde{W}_{j-1}$, which implies $T_{\tilde{W}_{j-1}} = \tilde{T}_{W_{j-1}} Tr_j$ in $V_j$.

Since both $\tilde{T}_{W_{j-1}}$ and $Tr_j$ in (D2) are invertible linear transformations, $T_{\tilde{W}_{j-1}}$ is an invertible linear transformation of $\tilde{W}_{j-1}$ onto $l^2(\mathbb{R})$. Since the translation $Tr_j$ is an isometry of $\tilde{W}_{j-1}$ onto $W_{j-1}$,

$$\{\psi_{1/2}(2^{j-1}x - k)\}_{k \in \mathbb{Z}} = \{\psi(2^{j-1}x - k + 1/2)\}_{k \in \mathbb{Z}} \tag{D3}$$

is a Riesz basis of $\tilde{W}_{j-1}$ where $\psi(x)$ is a wavelet of $W_0$. The space $\tilde{W}_{j-1}$ is obviously a reproducing kernel Hilbert space whose reproducing kernel has the form (19) with $\psi_{1/2}(x)$ in place of $\gamma(x)$, i.e. the kernel of $\tilde{W}_{j-1}$ can be written as $2^j \sum_{n \in \mathbb{Z}} \psi_{1/2}(2^{j-1}x - n) \tilde{\psi}_{1/2}(2^{j-1}x - n)$, where $\tilde{\psi}_{1/2}(2^{j-1}x)$ is the dual function of $\psi_{1/2}(2^{j-1}x)$ in $\tilde{W}_{j-1}$.

Since the data $T_{\tilde{W}_{j-1}}(f_s) = \{f_s(k/2^{j-1})\}_{k \in \mathbb{Z}}$ are the samples of $f_s(x)$ located at $x = k/2^{j-1}$ with $k \in \mathbb{Z}$, it follows from Proposition 2 that there exists an interpolation basis relative to $\{f_s(k/2^{j-1})\}_{k \in \mathbb{Z}}$ in $\tilde{W}_{j-1}$.

Let $\{S^{\psi}_{1/2}(2^{j-1}x - k)\}_{k \in \mathbb{Z}}$ be the interpolation basis relative to $\{f_s(k/2^{j-1})\}_{k \in \mathbb{Z}}$ in $\tilde{W}_{j-1}$. It follows from (20) that the inverse of $T_{\tilde{W}_{j-1}}$ is represented as

$$f_s(x) = \sum_{k \in \mathbb{Z}} f_s(k/2^{j-1}) S^{\psi}_{1/2}(2^{j-1}x - k) \tag{D4}$$

for every $f_s(x) \in \tilde{W}_{j-1}$.

Since $f_s(x + \frac{1}{2^j}) \in \tilde{W}_{j-1}$ for every element $f_s(x) \in W_{j-1}$, it follows from (D4) that

$$f_s(x + \frac{1}{2^j}) = \sum_k f_s(\frac{k}{2^{j-1}} + \frac{1}{2^j}) S^{\psi}_{1/2}(2^{j-1}x - k). \tag{D5}$$

The change of variable $x = y - 1/2^j$ in (D5) gives

$$f_s(y) = \sum_k f_s(\frac{k}{2^{j-1}} + \frac{1}{2^j}) S^{\psi}_{1/2}(2^{j-1}(y - 1/2^j) - k) = \sum_k f_s(\frac{k}{2^{j-1}} + \frac{1}{2^j}) S^{\psi}_{1/2}(2^{j-1}y - 1/2 - k). \tag{D6}$$

Since the operator $Tr_j$ is a bijective isometry of $\tilde{W}_{j-1}$ onto $W_{j-1}$ and since $S^{\psi}_{1/2}(2^{j-1}y)$ is a Riesz basis for $\tilde{W}_{j-1}$,

$$\{S^{\psi}(2^{j-1}x - k)\}_{k \in \mathbb{Z}} = Tr_j\{S^{\psi}_{1/2}(2^{j-1}x - k)\}_{k \in \mathbb{Z}} = \{S^{\psi}_{1/2}(2^{j-1}x - 1/2 - k)\}_{k \in \mathbb{Z}} \tag{D7}$$

is a Riesz basis for $W_{j-1}$.

By substituting (D7) into (D6), we have

$$f_s(x) = \sum_k f_s(\frac{k}{2^{j-1}} + \frac{1}{2^j}) S^{\psi}(2^{j-1}x - k), \tag{D8}$$

which proves (27). Equation (D8) implies that $\{S^{\psi}(2^{j-1}x - k)\}_k$ is just an interpolation wavelet basis relative to the samples of $f_s(x)$ located at $x = (2k+1)/2^j$.

Since $\psi(x)$ is an element of $W_0$, (D8) indicates that

$$\psi(x) = \sum_{k \in \mathbb{Z}} \psi(k + \frac{1}{2}) S^{\psi}(x - k). \tag{D9}$$

The Fourier transform of (D9) yields

$$\hat{S}^{\psi}(w) = \frac{\hat{\psi}(w)}{\sum_{k}\psi(k+1/2)e^{-iwk}}, \tag{D10}$$

which proves (28).

**Part 2: Proof that** $\sum_{k\in\mathbb{Z}}S^{\psi}(k/2)e^{-iwk/2}$ ( $\sum_{k\in\mathbb{Z}}\psi(k/2)e^{-iwk/2}$ ) **and** $\sum_{k\in\mathbb{Z}}\hat{S}^{\psi}(w+4k\pi)$ ( $\sum_{k\in\mathbb{Z}}\hat{\psi}(w+4k\pi)$ ) **converge pointwise everywhere.**

It follows from the first identity in (29) that

$$\frac{1}{2}\sum_{k\in\mathbb{Z}}\psi(k/2)e^{-iwk/2} = \frac{1}{2}\sum_{k\in\mathbb{Z}}\sum_{n\in\mathbb{Z}}q_n\phi(k-n)e^{-iwk/2} = Q_\psi(z)\sum_{k\in\mathbb{Z}}\phi(k)e^{-iwm/2}. \tag{D11}$$

Since we have assumed in (4) that $\sum_{k\in\mathbb{Z}}\phi(k)e^{-iwm/2}$ converges pointwise everywhere and since $Q_\psi(z)$ also converges pointwise everywhere, $Q_\psi(z)\sum_{k\in\mathbb{Z}}\phi(k)e^{-iwm/2}$ is a convergent series in (D11). Hence, equation (D11) indicates that $\sum_{k\in\mathbb{Z}}\psi(k/2)e^{-iwk/2}$ converges pointwise everywhere.

On the other hand, it follows from the first identity in (4) and the second identity in (29) that

$$Q_\psi(z)\sum_{k\in\mathbb{Z}}\phi(k)e^{-iwm/2} = Q_\psi(z)\sum_{k\in\mathbb{Z}}\hat{\phi}(w/2+2k\pi) = \sum_{k\in\mathbb{Z}}\hat{\psi}(w+4k\pi). \tag{D12}$$

Since $Q_\psi(z)\sum_{m\in\mathbb{Z}}\phi(k)e^{-iwm/2}$ converges pointwise everywhere, equation (D12) indicates that $\sum_{k\in\mathbb{Z}}\hat{\psi}(w+4k\pi)$ also does. By substituting (D12) into (D11), we have

$$\frac{1}{2}\sum_{k\in\mathbb{Z}}\psi(k/2)e^{-iwk/2} = \sum_{k\in\mathbb{Z}}\hat{\psi}(w+4k\pi), \tag{D13}$$

which proves the first identity in (30).

Since $\{\psi(k/2)\}_{k\in\mathbb{Z}} = \{\psi(k+1/2)\}_{k\in\mathbb{Z}} \cup \{\psi(k)\}_{k\in\mathbb{Z}}$ and since $\{e^{-iwk/2}\}_{k\in\mathbb{Z}}$ are linearly independent, it follows from (D13) that

$$\begin{cases} \sum_{k\in\mathbb{Z}}\psi(k)e^{-iwk} = \sum_{k\in\mathbb{Z}}\hat{\psi}(w+2k\pi) \\ e^{-iw/2}\sum_{k\in\mathbb{Z}}\psi(k+1/2)e^{-iwk} = \sum_{k\in\mathbb{Z}}(-1)^k\hat{\psi}(w+2k\pi), \end{cases} \tag{D14}$$

(also see (A2) and (A3) in Appendix A). Inserting the second identity in (D14) into (D10) yields (31).

Now, we prove the second identity in (30).

Since $\psi_{1/2}(x) = \psi(x+1/2)$, it follows from (D13) that

$$\begin{aligned}\sum_{k\in\mathbb{Z}}\psi_{1/2}(k/2)e^{-iwk/2} &= \sum_{k\in\mathbb{Z}}\psi(k/2+1/2)e^{-iwk/2} \\ &= e^{iw/2}\sum_{k\in\mathbb{Z}}\psi(k/2)e^{-iwk/2} = 2e^{iw/2}\sum_{k\in\mathbb{Z}}\hat{\psi}(w+4k\pi)\end{aligned}. \tag{D15}$$

Simultaneously, we have

$$\hat{\psi}_{1/2}(w+4k\pi) = e^{iw/2}\hat{\psi}(w+4k\pi) \tag{D16}$$

for every $k \in \mathbb{Z}$.

Inserting (D16) into (D15) yields

$$\frac{1}{2}\sum_{k\in\mathbb{Z}}\psi_{1/2}(k/2)e^{-iwk/2} = \sum_{k\in\mathbb{Z}}\hat{\psi}_{1/2}(w+4k\pi). \tag{D17}$$

Equation (D17) indicates that the Poisson summation formula holds for $\{\psi_{1/2}(k/2)\}_k$ when the Poisson summation formula holds for $\{\psi(k/2)\}_k$.

Since $\{S^{\psi}_{1/2}(x-k)\}_{k\in\mathbb{Z}}$ is an interpolation basis relative to $\{f_s(k)\}_{k\in\mathbb{Z}}$ in $\tilde{W}_0$ by (D4) and since $\{\psi_{1/2}(x-k)\}_k$ has formed a Riesz basis for $\tilde{W}_0$, it follows from Proposition 1 that

$$\frac{1}{2}\sum_{k\in\mathbb{Z}}S^{\psi}_{1/2}(k/2)e^{-iwk/2} = \sum_{k\in\mathbb{Z}}\hat{S}^{\psi}_{1/2}(w+4k\pi) \tag{D18}$$

if (D17) holds. It follows from (D7) that

$$\begin{cases} S^{\psi}_{1/2}(x) = S^{\psi}(x+1/2) \\ \hat{S}^{\psi}_{1/2}(w) = e^{iw/2}\hat{S}^{\psi}(w) \end{cases}. \tag{D19}$$

Inserting (D19) into (D18) yields the second identity in (30).

### Appendix E : PROOF OF LEMMA 2

**Proof.** By Parseval's identity,

$$\int_{-\infty}^{+\infty}\psi(x)\overline{\phi(x-n)}dx = \frac{1}{2\pi}\int_{-\infty}^{+\infty}\hat{\psi}(w)\overline{\hat{\phi}(w)}e^{iwn}dw = \frac{1}{2\pi}\sum_{k\in\mathbb{Z}}\int_{-\pi+2k\pi}^{+\pi+2k\pi}\hat{\psi}(w)\overline{\hat{\phi}(w)}e^{iwn}dw$$
$$= \frac{1}{2\pi}\int_{-\pi}^{\pi}\sum_{k\in\mathbb{Z}}\hat{\psi}(w+2k\pi)\overline{\hat{\phi}(w+2k\pi)}e^{iwn}dw \qquad, \text{ for } n\in\mathbb{Z}. \tag{E1}$$

Since $\psi(x)$ is an element of $W_0$ and since $W_0$ is orthogonal to $V_0$, we have

$$\int_{-\infty}^{+\infty}\psi(x)\overline{\phi(x-n)}dx = 0. \tag{E2}$$

Inserting (E2) into (E1) yields

$$\int_{-\pi}^{\pi}\sum_{k\in\mathbb{Z}}\hat{\psi}(w+2k\pi)\overline{\hat{\phi}(w+2k\pi)}e^{iwn}dw = 0. \tag{E3}$$

Since $\{e^{iwn}\}_n$ is a complete orthonormal basis of $L^2[-\pi,\pi]$, it follows from (E3) that

$$\sum_{k\in\mathbb{Z}}\hat{\psi}(w+2k\pi)\overline{\hat{\phi}(w+2k\pi)} = 0 \tag{E4}$$

almost everywhere.

It follows from (8) and (38) that

$$\begin{cases} \hat{\phi}(w+4k\pi) = P_\phi\left(e^{-\frac{i(w+4k\pi)}{2}}\right)\hat{\phi}(\frac{w}{2}+2k\pi) \\ \qquad\qquad = P_\phi(z)\hat{\phi}(\frac{w}{2}+2k\pi), \quad k \in \mathbb{Z} \\ \hat{\phi}(w+2\pi+4k\pi) = P_\phi\left(e^{-\frac{i(w+4k\pi+2\pi)}{2}}\right)\hat{\phi}(\frac{w}{2}+\pi+2k\pi) \\ \qquad\qquad = P_\phi(-z)\hat{\phi}(\frac{w}{2}+2k\pi+\pi), \quad k \in \mathbb{Z} \end{cases} \quad (E5)$$

and

$$\begin{cases} \hat{\psi}(w+4k\pi) = Q_\psi(z)\hat{\phi}(\frac{w}{2}+2k\pi), \quad k \in \mathbb{Z} \\ \hat{\psi}(w+2\pi+4k\pi) = Q_\psi(-z)\hat{\phi}(\frac{w}{2}+2k\pi+\pi), \quad k \in \mathbb{Z} \end{cases} \quad (E6)$$

Equations (E5) and (E6) indicate that

$$\sum_{k\in\mathbb{Z}} \hat{\psi}(w+2k\pi)\overline{\hat{\phi}(w+2k\pi)} = \sum_{k\in\mathbb{Z}} \hat{\psi}(w+4k\pi)\overline{\hat{\phi}(w+4k\pi)} + \sum_{k\in\mathbb{Z}} \hat{\psi}(w+2\pi+4k\pi)\overline{\hat{\phi}(w+2\pi+4k\pi)}$$
$$= Q_\psi(z)\sum_{k\in\mathbb{Z}}\left|\hat{\phi}(\frac{w}{2}+2k\pi)\right|^2 \overline{P_\phi(z)} + Q_\psi(-z)\sum_{k\in\mathbb{Z}}\left|\hat{\phi}(\frac{w}{2}+2k\pi+\pi)\right|^2 \overline{P_\phi(-z)} \quad . \quad (E7)$$

Inserting (41) and (E4) into (E7) yields

$$Q_\psi(z)E_\phi(z)\overline{P_\phi(z)} + Q_\psi(-z)E_\phi(e^{-i(w+2\pi)/2})\overline{P_\phi(-z)} = Q_\psi(z)E_\phi(z)\overline{P_\phi(z)} + Q_\psi(-z)E_\phi(-z)\overline{P_\phi(-z)} = 0,$$

which verifies (40).□

### Appendix F : PROOF OF THEOREM 2

since $\{S^\phi(2^j x-k)\}_{k\in\mathbb{Z}}$ and $\{S^\psi(2^j x-k)\}_{k\in\mathbb{Z}}$ respectively form Riesz bases for $V_j$ and $W_j$, $P_s(z)$ in (12) and $Q_s(z)$ in (34) form a pair of reconstruction filters in an MRA $\{V_j\}_{j\in\mathbb{Z}}$.

**First, we show that $P_s(-z) + P_s(z) = 1$ and $Q_s(z)/z - Q_s(-z)/z = 1$.**

It follows from (12) and (34) that

$$\begin{cases} P_s(z) = \sum_{k\in\mathbb{Z}} \hat{S}^\phi(w+4k\pi) \\ Q_s(z)/z = \frac{1}{z}\sum_{k\in\mathbb{Z}} \hat{S}^\psi(w+4k\pi) \end{cases} . \quad (F1)$$

The change of variable $w = u + 2\pi$ in (F1) gives

$$\begin{cases} P_s(-z) = \sum_{k\in\mathbb{Z}} \hat{S}^\phi(u+\pi+4k\pi) \\ Q_s(-z)/z = \frac{1}{z}\sum_{k\in\mathbb{Z}} \hat{S}^\psi(u+\pi+4k\pi) \end{cases} . \quad (F2)$$

It follows from (F1) and (F2) that

$$P_s(z) + P_s(-z) = \sum_{k \in \mathbb{Z}} \hat{S}^\phi(w + 4k\pi) + \sum_{k \in \mathbb{Z}} \hat{S}^\phi(u + 2\pi + 4k\pi) = \sum_{k \in \mathbb{Z}} \hat{S}^\phi(u + 2\pi k) \tag{F3}$$

and

$$Q_s(z)/z - Q_s(-z)/z = \frac{1}{z}(\sum_{k \in \mathbb{Z}} \hat{S}^\psi(w + 4k\pi) - \sum_{k \in \mathbb{Z}} \hat{S}^\psi(u + 2\pi + 4k\pi))$$

$$= \frac{1}{z} \sum_{k \in \mathbb{Z}} (-1)^k \hat{S}^\psi(w + 2k\pi) = \sum_{k \in \mathbb{Z}} \hat{S}^\psi(w + 2k\pi) e^{i(w + 2k\pi)/2} . \tag{F4}$$

Applying (10) and (31) respectively to (F3) and (F4) yields

$$\begin{cases} P_s(z) + P_s(-z) = \sum_{k \in \mathbb{Z}} \hat{S}^\phi(u + 2\pi k) = 1 \\ Q_s(z)/z - Q_s(-z)/z = \sum_{k \in \mathbb{Z}} \hat{S}^\psi(w + 2k\pi) e^{i(w + 2k\pi)/2} = 1 . \end{cases} \tag{F5}$$

It follows from the second identity in (F5) that $Q_s(1) - Q_s(-1) = 1$. Recall that as a pair of reconstruction filter, $(P_s, Q_s)$ satisfy the conditions

$$\begin{cases} P_s(1) = 1 \quad P_s(-1) = 0 \\ Q_s(1) = 0 \end{cases} \tag{F6}$$

(Equation 5.4.1 in [20]). Hence, equations (F5) and (F6) imply (43).

**Second, we show that there exist two constants $A_s$ and $B_s$ with $0 < A_s \leq B_s < +\infty$ such that**

$$A_s < |PE_s(w)|^2 = |\overline{P_s(-z)}E_s(-z) + \overline{P_s(z)}E_s(z)|^2 < B_s .$$

It follows from the second identity in (F5) that

$$Q_s(z)\left(\overline{P_s(-z)}E_s(-z) + E_s(z)\overline{P_s(z)}\right) = Q_s(z)\overline{P_s(-z)}E_s(-z) + Q_s(z)E_s(z)\overline{P_s(z)}$$

$$= z\overline{P_s(-z)}E_s(-z) + Q_s(-z)\overline{P_s(-z)}E_s(-z) + Q_s(z)E_s(z)\overline{P_s(z)} \tag{F7}$$

Since $P_s(z)$ in (12) and $Q_s(z)$ in (34) form a pair of reconstruction filters in an MRA $\{V_j\}_{j \in \mathbb{Z}}$, it follows from (40) that

$$Q_s(z)E_s(z)\overline{P_s(z)} + Q_s(-z)E_s(-z)\overline{P_s(-z)} = 0 . \tag{F8}$$

Inserting (F8) into (F7) yields

$$Q_s(z)\left(\overline{P_s(-z)}E_s(-z) + E_s(z)\overline{P_s(z)}\right) = z\overline{P_s(-z)}E_s(-z) \tag{F9}$$

On the other hand, applying (F5) to (37) yields

$$\Delta_{P_s, Q_s} = P_s(z)Q_s(z) - P_s(-z)Q_s(z) - zP_s(z) . \tag{F10}$$

Now, consider the following formula

$$\Delta_{P_s, Q_s}(E_s(-z)\overline{P_s(-z)} + E_s(z)\overline{P_s(z)}) . \tag{F11}$$

Inserting (F10) into (F11) yields

$$\Delta_{P_s,Q_s}\left(E_s(-z)\overline{P_s(-z)}+E_s(z)\overline{P_s(z)}\right)=\left(P_s(z)-P_s(-z)\right)Q_s(z)\left(E_s(-z)\overline{P_s(-z)}+E_s(z)\overline{P_s(z)}\right)$$
$$-zP_s(z)\left(E_s(-z)\overline{P_s(-z)}+E_s(z)\overline{P_s(z)}\right) \quad . \tag{F12}$$

Applying (F9) to (F12) yields

$$\Delta_{P_s,Q_s}\left(E_s(-z)\overline{P_s(-z)}+E_s(z)\overline{P_s(z)}\right)=\left(P_s(z)-P_s(-z)\right)z\overline{P_s(-z)}E_s(-z)$$
$$-zP_s(z)\left(E_s(-z)\overline{P_s(-z)}+E_s(z)\overline{P_s(z)}\right)=-z\left(E_s(-z)|P_s(-z)|^2+E_s(z)|P_s(z)|^2\right). \tag{F13}$$

It follows from (16) and (13) that

$$E_s(-z)|P_s(-z)|^2+E_s(z)|P_s(z)|^2 = \sum_{k=-\infty}^{+\infty}\left|\hat{S}^\phi\left(\frac{w}{2}+2k\pi+\pi\right)P_s(-z)\right|^2 + \sum_{k=-\infty}^{+\infty}\left|\hat{S}^\phi\left(\frac{w}{2}+2k\pi\right)P_s(z)\right|^2$$
$$= \sum_{k=-\infty}^{+\infty}\left|\hat{S}^\phi(w+2k\pi)\right|^2 = E_s(z^2) \quad . \tag{F14}$$

Inserting (F14) into (F13), we have

$$-zE_s(z^2) = \Delta_{P_s,Q_s}\left(E_s(-z)\overline{P_s(-z)}+E_s(z)\overline{P_s(z)}\right). \tag{F15}$$

Since $(P_s(z),Q_s(z))$ forms a pair of reconstruction filters, equation (36) holds. Hence, it follows from (F15) that

$$E_s(-z)\overline{P_s(-z)}+E_s(z)\overline{P_s(z)} = -zE_s(z^2)/\Delta_{P_s,Q_s} \tag{F16}$$

Since the sequence $\{\phi(x-k)\}_k$ is a frame, there should exist two constants $A_{Es}$ and $B_{Es}$ with $0 < A_{Es} \le B_{Es} < +\infty$ such that

$$A_{Es} \le E(z^2) \le B_{Es} \text{ (Theorem 9 in [19])}. \tag{F17}$$

Applying (36) and (F17) to (F16) yields

$$A_{Es}^2/B_\Delta \le \left|E_s(-z)\overline{P_s(-z)}+E_s(z)\overline{P_s(z)}\right|^2 = \left|E_s(z^2)/\Delta_{P_s,Q_s}\right|^2 \le B_{Es}^2/A_\Delta, \tag{F18}$$

which proves (23).

Equation (F18) implies that

$$\frac{1}{\overline{P_s(-z)}E_s(-z)+E_s(z)\overline{P_s(z)}} \in L^2[-\pi,\pi]. \tag{F19}$$

Hence, it follows from (F9) and (F19) that

$$Q_s(z) = \frac{zE_s(-z)\overline{P_s(-z)}}{\overline{P_s(-z)}E_s(-z)+E_s(z)\overline{P_s(z)}}, \tag{F20}$$

which verifies equation (44).□

## Appendix G  A LIST OF SYMBOLS

| Symbol | Description | Symbol | Description |
|---|---|---|---|
| $W_j$ | Wavelet space | $\hat{S}^{\psi}_{1/2}(w)$ | The Fourier transform of $S^{\psi}_{1/2}(x)$ |
| $V_j$ | Approximation space of multiresolution analysis | $\{\gamma(t-k)\}_{k\in\mathbb{Z}}$ | Riesz basis of a Hilbert space |
| $U_j(\gamma)$ | Space spanned by a Riesz basis $\{\gamma(2^j x - k)\}_k$ | $\{\tilde{\gamma}(t-k)\}_{k\in\mathbb{Z}}$ | Dual basis of $\{\gamma(t-k)\}_{k\in\mathbb{Z}}$ |
| $\tilde{W}_{j-1}$ | $\tilde{W}_{j-1} = \{f(x)\mid f(x-1/2^j)\in W_{j-1}\}$ | $\hat{\gamma}(w)$ | The Fourier transform of $\gamma(t)$ |
| $L^2(\mathbb{R})$ | Square integrable function | $\{S^{\phi}(2^j x - k)\}_{k\in\mathbb{Z}}$ | Interpolation basis relative to the samples $\{f_s(k/2^j)\}_{k\in\mathbb{Z}}$ for $V_j$ |
| $l^2(\mathbb{R})$ | Finite energy discrete signals $\sum_{k=-\infty}^{+\infty}|f_s(k)|^2 < +\infty$ | $\{\tilde{S}^{\phi}(2^j x - k)\}_{k\in\mathbb{Z}}$ | Dual basis of $\{S^{\phi}(2^j x - k)\}_{k\in\mathbb{Z}}$ for $V_j$ |
| $f_s(x)$ | Signal to be recovered | $\{S^{\psi}(2^j x - k)\}_{k\in\mathbb{Z}}$ | Interpolation basis relative to the samples $\{f_s(k/2^j+1/2^{j+1})\}_{k\in\mathbb{Z}}$ for $W_j$ |
| $f_{ap}(x)$ | Recovery of signal | $\{\tilde{S}^{\psi}(2^j x - k)\}_{k\in\mathbb{Z}}$ | Dual basis of $\{S^{\psi}(2^j x - k)\}_{k\in\mathbb{Z}}$ for $W_j$ |
| $\hat{f}_s(w)$ | Fourier transform of $f_s(x)$ | $q_j(x,y)$ | Reproducing kernel defined in (19) |
| $\phi(x)$ | Scaling function of $V_0$ | $P_{\phi}(z)$ | $l^2$-sequence defined in (9) |
| $\hat{\phi}(w)$ | The Fourier transform of $\phi(x)$ | $P_s(z)$ | $l^2$-sequence defined in (12) |
| $\psi(x)$ | Wavelet of $W_0$ | $\tilde{Q}_s(z)$ | $l^2$-sequence defined in (C25) |
| $\hat{\psi}(w)$ | The Fourier transform of $\psi(x)$ | $Q_{\psi}(z)$ | $l^2$-sequence defined in (39) |
| $S^{\phi}(x)$ | Interpolation function for $V_0$ | $Q_s(z)$ | $l^2$-sequence defined in (34) |
| $\hat{S}^{\phi}(w)$ | The Fourier transform of $S^{\phi}(x)$ | $E_s(z)$ | Period function defined in (16) |
| $S^{\psi}(x)$ | Interpolation function for $W_0$ | $E_{\phi}(z)$ | Period function defined in (41) |
| $\hat{S}^{\psi}(w)$ | The Fourier transform of $S^{\psi}(x)$ | $\Delta_{P_s,Q_s}$ | Function defined in (37) |
| $\overline{f}(x)$ | Complex conjugate of $f(x)$ | $T_{\tilde{W}_{j-1}}$ | Sampling operator $T_{\tilde{W}_{j-1}}(f_s) = \{f_s(k/2^{j-1})\}_{k\in\mathbb{Z}}$ from $\tilde{W}_{j-1}$ to $l^2(\mathbb{R})$ |
| $S^{\phi}_{j,k}(x)$ | $S^{\phi}(2^j x - k)$ | $T_{V_0}$ | Sampling operator $T_{V_0}(f_s) = \{f_s(k)\}_k$ from $V_0$ to $l^2(\mathbb{R})$ |
| $\psi_Q(x)$ | Wavelet defined in (C26) | $T_{U_j}$ | Sampling operator $T_{U_j}(f_s) = \{f_s(k/2^j)\}_k$ from $U_j$ to $l^2(\mathbb{R})$ |
| $\hat{\psi}_Q(w)$ | The Fourier transform of $\psi_Q(x)$ | $\tilde{T}_{W_{j-1}}$ | Sampling operator $\tilde{T}_{W_j}(f_s) = \{f_s(k/2^j+1/2^{j+1})\}_{k\in\mathbb{Z}}$ from $W_j$ to $l^2(\mathbb{R})$ |
| $\psi_{1/2}(x)$ | Function defined in (D3) | $Tr_j$ | Shift operator defined in (D1) |
| $\hat{\psi}_{1/2}(w)$ | The Fourier transform of $\psi_{1/2}(x)$ | $\oplus$ | Orthogonal sum. |
| $S^{\psi}_{1/2}(x)$ | Interpolation function of $\tilde{W}_0$ defined in (D4) | | |